\documentclass[a4paper,11pt]{article}
\usepackage{setspace}
\usepackage[margin=2.5cm]{geometry}
\usepackage[table]{xcolor}
\usepackage{colortbl}
\usepackage[ruled,vlined]{algorithm2e}
\usepackage{amsmath}
\usepackage{amsthm}
\usepackage{amssymb}
\usepackage{graphicx}
\usepackage{placeins}
\usepackage{longtable}
\usepackage{comment}
\PassOptionsToPackage{hyphens}{url}
\usepackage{hyperref}
\usepackage{booktabs}
\usepackage{array}
\usepackage{microtype}
\usepackage{caption}
\usepackage{multirow}
\usepackage[authoryear,round]{natbib}
\usepackage{authblk}

\title{Constraint-Relaxation Trade-offs as a By-Product of \\Evolutionary Black-Box Optimization}
\author[a]{Giovanni Misitano\thanks{Corresponding author. E-mail: \texttt{giovanni.a.misitano@jyu.fi}}}
\author[b]{Juergen Branke}
\affil[a]{University of Jyv\"askyl\"a, Faculty of Information Technology, P.O. Box 35,
      FI-40014 University of Jyv\"askyl\"a, Finland}
\affil[b]{Warwick Business School, The University of Warwick, Coventry, CV4 7AL,
      United Kingdom}
\date{\today}

\begin{document}

\maketitle
\begin{abstract}
      Constraints in practical optimization problems are not always fixed:\ a
      budget can be stretched, a tolerance loosened. A decision maker may
      therefore want to know how much the objective would improve if a
      constraint were relaxed. Linear and smooth nonlinear programming answer
      this with the shadow price, the constraint's dual variable. For
      black-box problems no dual information is available, and the gain is often a
      nonlinear function of the amount of relaxation. We therefore target the
      shadow-front:\ the trade-off between how far a constraint is relaxed, up
      to a tolerance the decision maker accepts, and how much the objective
      can improve. Four scenarios arise, depending on whether the end points
      or the whole curve matter and whether the constraints are relaxed one at
      a time or jointly. We show that this information can be obtained as a
      free by-product of a standard feasibility-first evolutionary algorithm by modifying only
      its environmental selection operator. We propose four such operators,
      the relaxed, joint, ranking, and hypervolume modes, and compare them
      with a feasibility-first baseline on ten benchmark problems. As expected, all modes
      recover the shadow-front better than the baseline in most configurations. Perhaps surprisingly,
      they also match or outperform the baseline on the original constrained
      problem in most configurations, and so obtain the information essentially
      for free. The hypervolume mode
      is, overall,
      the most accurate in all four scenarios, at a cost that grows steeply
      with the number of constraints; the relaxed, joint, and ranking modes
      are cheap alternatives. \end{abstract}


\section{Introduction}
\label{sec:introduction}

Most optimization problems are solved as if their constraints were fixed. In
practice, many constraints are negotiable:\ a production limit can be raised by
hiring staff or working overtime, a design tolerance can be loosened at some
cost, a budget can be stretched. A decision maker then may wish to know more than
the constrained optimum. They want to know what a relaxation would buy:\ how
much the objective could improve if a given constraint were relaxed by a given
amount. In linear and smooth nonlinear programming this question has a
classical answer. The shadow price of a constraint, its optimal dual variable or
Lagrange multiplier, gives the marginal improvement of the optimal objective
value per unit of relaxation~\citep{bertsekas1999nonlinear,boyd2004convex}, and
every linear programming solver reports it at no extra
cost. Shadow prices are central to post-optimality analysis and are widely used
in cost-benefit analysis, pricing, and resource allocation~\citep{mas1995microeconomic}.

However, many problems in engineering and operations are black boxes:\ the
objective and the constraints are evaluated by a simulator or an experiment,
and no gradients, convexity, or closed forms exist. Such problems are
commonly solved with metaheuristics such as evolutionary algorithms, which return a feasible
near-optimal solution. Dual variables are not readily available, and even the notion of a single
marginal rate is inadequate, because the attainable gain is, in general, a
nonlinear function of the amount of relaxation. The obvious remedies are
unattractive:\ re-optimizing for each relaxation level multiplies the evaluation
cost, and reformulating the constraints as objectives and solving the resulting
multiobjective problem spends the budget on the entire trade-off surface, far
beyond the small relaxations a decision maker would realistically accept.

This paper proposes a different approach. We replace the infinitesimal relaxation
of the classical shadow price by a finite one:\ for each constraint the decision
maker specifies the largest relaxation they are willing to consider, and we
seek the \emph{shadow-front}, the non-dominated trade-off between the amount of
relaxation, up to that tolerance, and the attainable objective value. What a
decision maker may request varies along two dimensions. First,
they may want only its end points, that is, the gain obtained at the full
relaxation; or the whole curve between them, to judge where the relaxation
stops paying off. As an example, a plant manager asking what an extra shift on a
line would buy wants a single number; a designer deciding how far to loosen an
allowable stress wants to see how the weight changes along the way. Second, they
may relax one constraint at a time, because each has its own owner or cost, or
all of them jointly when, as in a negotiation, all limits can be traded off against each other. The
two dimensions give four scenarios, each demanding a different emphasis in
the search.

We show that in each scenario the required information can be obtained as a
free by-product of a standard evolutionary algorithm, by modifying only
its environmental selection operator and without any additional objective
evaluations. The evaluations that any constraint-handling scheme spends on
solutions near the feasibility boundary are steered, and recorded, so
that they map out the shadow-front. Our contributions are as follows:
\begin{enumerate}
      \item We adapt the concept of a shadow price to constrained
      single-objective black-box optimization through the relaxation gain, the
      objective improvement obtained from a finite relaxation, and the
      shadow-front, and we structure the information a decision maker may want
      from them into four scenarios (Section~\ref{subsec:definitions}).
      
      \item We propose four environmental selection operators that recover this
      information during a single run:\ a \emph{relaxed} mode that interleaves
      per-constraint relaxed rankings and targets the end points under single
	relaxations; a \emph{joint} mode similar to the relaxed mode, but
	with a ranking that considers all constraints relaxed simultaneously and targets the end point
	at the joint relaxations;
      a \emph{ranking} mode that treats these rankings of the relaxed mode as
      objectives of a multiobjective problem and covers the front under single
      relaxations; and a \emph{hypervolume} mode that selects on contribution to
      the hypervolume of the shadow-front under joint relaxation
      (Section~\ref{sec:method}).
      \item We evaluate the four modes against a feasibility-first
      baseline, which stands in for conventional constraint handling, on ten
      benchmark problems, three population sizes, and three relaxation levels,
      measuring both optimization performance and the accuracy of the
      recovered shadow-front
      (Section~\ref{sec:experiments}).
\end{enumerate}
All four modes recover the shadow-front better than the baseline, and on most
problems they match or improve on it in finding the constrained optimum. Overall, the hypervolume mode is the most
accurate in all four scenarios, both for the end points and for the trade-off curves,
but its selection cost grows steeply with the number of constraints and the
population size, and with three or four constraints and the largest population it is by far the most expensive. Of
the cheaper modes, the relaxed mode is the better
estimator of the end points when constraints are relaxed individually, and the
joint mode the best estimator of the end point when all constraints are relaxed
simultaneously. The ranking mode is the strongest of the cheaper modes,
covering the trade-off region best.
To our knowledge, no prior method reports the objective gain of relaxing
a given constraint by a specified amount as an output of a single black-box
optimization run. The method,
experimental pipeline, and data are made openly available and are fully
reproducible.

The remainder of the paper is organized as follows. Section~\ref{sec:background}
reviews the concepts we build on and the related literature.
Section~\ref{sec:method} defines the relaxation gain and the shadow-front for
black-box problems, lays out the four scenarios, and presents the four selection
modes. Section~\ref{sec:experiments} describes the experimental setup and
reports the results, which Section~\ref{sec:discussion} discusses, including
which mode to use in which scenario. Section~\ref{sec:conclusions} concludes.

\section{Background and related work}
\label{sec:background}

\subsection{Constrained black-box optimization and the shadow price}
\label{subsec:problem_and_shadow}

We consider constrained optimization problems with $k \geq 1$ inequality
constraints:
\begin{equation}
      \label{eq:problem_definition}
      \begin{split}
            \min_{\mathbf{x}} & \quad f(\mathbf{x})                                                   \\
            \text{s.t.}       & \quad c_j(\mathbf{x}) \le 0\quad \text{for}\; j=1,\dots,k             \\
                              & \quad \mathbf{x}^\text{lb} \leq \mathbf{x} \leq \mathbf{x}^\text{ub},
      \end{split}
\end{equation}
where $f:\mathbb{R}^n \to \mathbb{R}$ is the objective function,
$\mathbf{x} = [x_1,\dots,x_n] \in \mathbb{R}^n$ is the decision variable vector,
$c_j:\mathbb{R}^n \to \mathbb{R}$ are the constraint functions, and
$\mathbf{x}^\text{lb}, \mathbf{x}^\text{ub}$ are element-wise box-bounds. A solution
is feasible when all constraints hold. We assume $f$ and $c_j$ to be black
boxes, as in simulation- or experiment-based problems, so that gradients,
convexity, and closed analytic forms are unavailable.

The shadow price of a constraint is the marginal value of relaxing it:\ the
optimal dual variable, which gives the rate at which the optimal objective value
improves under an infinitesimal relaxation of the constraint's right-hand
side~\citep{bertsekas1999nonlinear,boyd2004convex}. Linear programming solvers
return shadow prices with the optimum, and in nonlinear programming they are the
optimal Lagrange multipliers under standard regularity conditions. For
evolutionary algorithms on smooth problems, the Karush--Kuhn--Tucker (KKT) proximity measure
of~\citet{tulshyan2010investigating,dutta2013approximate} yields approximate
Lagrange multipliers of near-optimal solutions as a by-product, but it needs the
gradients of the objective and the constraints. In the
black-box setting~\eqref{eq:problem_definition} no multipliers are available,
so any analogue must be constructed from evaluated solutions. This is
discussed further in Section~\ref{subsec:definitions}.

\subsection{Evolutionary algorithms and constraint handling}
\label{subsec:ea_and_constraints}

Evolutionary algorithms apply to a wide range of problems without much
problem-specific tailoring, and deliver good, though not necessarily optimal,
solutions within acceptable time~\citep{eiben2015introduction}. A generational evolutionary algorithm
improves a population of candidate solutions over generations. In each
generation, mating selection chooses parents, recombination (crossover) and
mutation create offspring from them, and environmental selection decides which
solutions from the combined pool of parents and offspring survive into the next
generation. Our contribution lies entirely in environmental selection:\ the
modes proposed in Section~\ref{sec:method} change only how the combined pool is
ranked and truncated, and leave every other operator untouched.

Constraints are typically also handled at this step. A standard approach, which
needs no penalty parameter~\citep{coello2002theoretical}, is the feasibility-first
rule of \citet{deb2000efficient}:\ any feasible
solution is preferred to any infeasible one, feasible solutions are ordered by
objective value, and infeasible solutions by their constraint violation. The
early survey in \citet{michalewicz1996evolutionary} groups constraint handling
into feasibility-preserving operators, penalty functions, methods separating
feasible from infeasible solutions, and hybrids; feasibility rules belong to the
third category and are documented as simple and widely
adopted~\citep{mezuramontes2011constraint,rahimi2023review}.
Alternatives balance objective value against violation stochastically
\citep{runarsson2000stochastic,runarsson2005search} or compare feasibility up to
a tolerance $\varepsilon$ that is annealed to zero during the
run~\citep{takahama2006constrained,takahama2010constrained}. Recent methods
adaptively anneal a violation boundary~\citep{chen2024acrea} or learn which
constraint-handling technique to apply~\citep{wang2024drl}; see
\citet{liang2024review} for a survey. All of these relax constraints only
temporarily, as a device to reach the feasible optimum, and none reports the
objective value of a relaxation.

Of these constraint-handling techniques, the $\varepsilon$-constrained
method~\citep{takahama2006constrained,takahama2010constrained} is mechanically
the closest to our relaxation thresholds. The difference is one of purpose:\ $\varepsilon$ is an
internal tolerance that vanishes by the end of the run, whereas our threshold is
a decision-maker-specified relaxation whose objective consequence is the output.

\subsection{Multiobjective optimization concepts}
\label{subsec:moo}

Two of our modes and our evaluation borrow from multiobjective
optimization~\citep{branke2008multiobjective}. In a multiobjective problem,
$m$ conflicting objectives $f_1,\dots,f_m$ are minimized simultaneously. A
solution $\mathbf{x}^1$ dominates $\mathbf{x}^2$ if $f_i(\mathbf{x}^1) \leq
      f_i(\mathbf{x}^2)$ for all $i$ and $f_j(\mathbf{x}^1) < f_j(\mathbf{x}^2)$ for
some $j$. The non-dominated solutions of a set are those not dominated by any
other member of the set; the non-dominated solutions of the whole feasible set
form the Pareto front. The shadow-front we seek is itself a non-dominated set,
with the amount of relaxation of each constraint and the objective value as its
axes, which is why the following tools apply to it.

Non-dominated sorting~\citep{deb2002fast} partitions a population into ranked
fronts:\ the non-dominated solutions form the first front, the non-dominated
solutions of the remainder the second, and so on. Within a front, a niching
measure favors solutions in sparsely populated regions; we use the mean
distance to the $\kappa$ nearest neighbors, a variant of the density estimate of
SPEA2~\citep{zitzler2001spea2}. The
hypervolume~\citep{zitzler1999multiobjective} measures the volume of objective
space that a set dominates relative to a reference point, typically a point slightly worse than the nadir
(the worst value of each objective over the Pareto front). It rewards both
convergence and spread, and it is the basis of indicator-based algorithms such
as SMS-EMOA~\citep{beume2007sms}, which select on each solution's contribution
to the hypervolume. We use it both to measure how well a mode covers the
shadow-front and, in the hypervolume mode, to drive selection.

\subsection{Related literature}
\label{subsec:related_literature}

Our modes (Section~\ref{sec:method}) are related to
\emph{multiobjectivization}, the reformulation of a
single-objective problem as a multiobjective one. Adding or decomposing
objectives can open monotonically improving paths and reduce local
optima~\citep{knowles2001multiobjectivization,jensen2004helper,handl2008decomposition}.
Applied to constraints, this gives methods that treat constraint violation as
one or more additional objectives:\ ranking feasible above infeasible solutions
by non-dominance~\citep{coello2000constrainthandling}, bi-objective formulations
of objective versus aggregated violation that keep infeasible solutions
alongside feasible ones if their objective is better, with infeasible-solution
archiving~\citep{cai2006multiobjective,wang2012combining}, self-adaptive
selection strategies~\citep{jiao2013novel}, and the early COMOGA
method~\citep{surry1997comoga}; see \citet{segura2016using} for a survey. 

Several constraint-handling approaches deliberately keep infeasible solutions
in the population~\citep[e.g.,][]{hamida2002aschea,mezura2005simple,isaacs2008blessings}.
Most notably, the infeasibility-driven evolutionary algorithm of
\citet{ray2009infeasibility} retains a fixed proportion of marginally
infeasible solutions, arguing that optima of constrained problems often lie on
the constraint boundary and that approaching it from both the feasible and the
infeasible side speeds up convergence. The authors also point out that such
marginally infeasible solutions may be valuable for trade-off studies, and
list such solutions for one problem with their objective values and
violations, but they neither bound the relaxation per constraint nor derive
the gain as a function of the relaxation. Their empirical finding that retaining
infeasible solutions improves rather than hinders the search for the
constrained optimum is consistent with our observation that shadow-front
information can be obtained essentially for free.

In all of the above cited papers, the aim is to reach the constrained
optimum; none bounds the relaxation per constraint by a user-specified threshold
or reports the resulting relaxation-versus-gain front. Several also aggregate
all violations into one objective, whereas we keep each constraint separate.
\citet{deb2010fast} argue that a penalty parameter whose
penalized minimum lies arbitrarily close to the constrained optimum equals the
constraint's Lagrange multiplier, i.e., its shadow price. In their hybrid method \citep{deb2010fast,deb2013biobjective} they fit a cubic curve to the non-dominated front of objective
versus the sum of the normalized constraint violations near zero violation, and
use the negative of its slope at zero violation to set the penalty parameter of a subsequent local search. 
Within constrained black-box optimization, \citet{ledigabel2024taxonomy}
classifies constraints; ours are quantifiable and relaxable in their terms.

The term shadow price itself has appeared in evolutionary and metaheuristic
optimization for other purposes. The shadow-price-guided genetic algorithm of
\citet{shen2011shadow,shen2011energy} attaches a potential-improvement score to
solution components to steer the operators. Closest to our goal is the genetic-algorithm
post-optimality analysis of \citet{saha2009ga}, which reports how the optimal
cost changes as constraint limits are relaxed, but obtains this by re-optimizing
at each limit. To our knowledge, no prior work delivers the objective
gain of relaxing a given constraint by a specified amount as a by-product
of a single black-box optimization run.

\section{Method}
\label{sec:method}

\subsection{Relaxation gain, shadow-front, and four scenarios}
\label{subsec:definitions}

For each constraint we introduce a relaxation threshold $t_j > 0$, the largest
relaxation of $c_j$ the decision maker is willing to consider, and collect them
in $\mathbf{t} = [t_1, \dots, t_k]$. In practice a threshold of zero is also
allowed, for a constraint the decision maker will not relax:\ that constraint is
then simply kept strict, as $c_4$ of the pressure vessel
(Table~\ref{tab:problems_thresholds}). For a relaxation vector $\mathbf{s}$ with
$0 \leq s_j \leq t_j$, let
\begin{equation}
      \label{eq:relaxed_optimum}
      f^\star(\mathbf{s}) = \min \{ f(\mathbf{x}) : c_j(\mathbf{x}) \le s_j \text{ for } j=1,\dots,k,\;
      \mathbf{x}^\text{lb} \leq \mathbf{x} \leq \mathbf{x}^\text{ub} \}
\end{equation}
be the optimal value under relaxation $\mathbf{s}$, so that $f^\star(\mathbf{0})$
is the optimum of the original problem. We define the \emph{relaxation gain} of
$\mathbf{s}$ as the resulting improvement,
\begin{equation}
      \label{eq:shadow_price}
      \pi(\mathbf{s}) = f^\star(\mathbf{0}) - f^\star(\mathbf{s}) \geq 0.
\end{equation}
The relaxation gain is the finite-relaxation counterpart of the shadow price.
Unlike the shadow price, it is a total gain, not a per-unit rate; in a
black-box problem the latter is neither available nor sufficient, since the
gain is in general nonlinear in $\mathbf{s}$. We reserve the term shadow price
for the classical dual quantity. A solution is \emph{threshold-feasible} if
$c_j(\mathbf{x}) \leq t_j$ for all $j$, and strictly feasible if it is feasible
for problem~\eqref{eq:problem_definition}.

\begin{figure}[htb]
      \centering
      \includegraphics[width=0.7\linewidth]{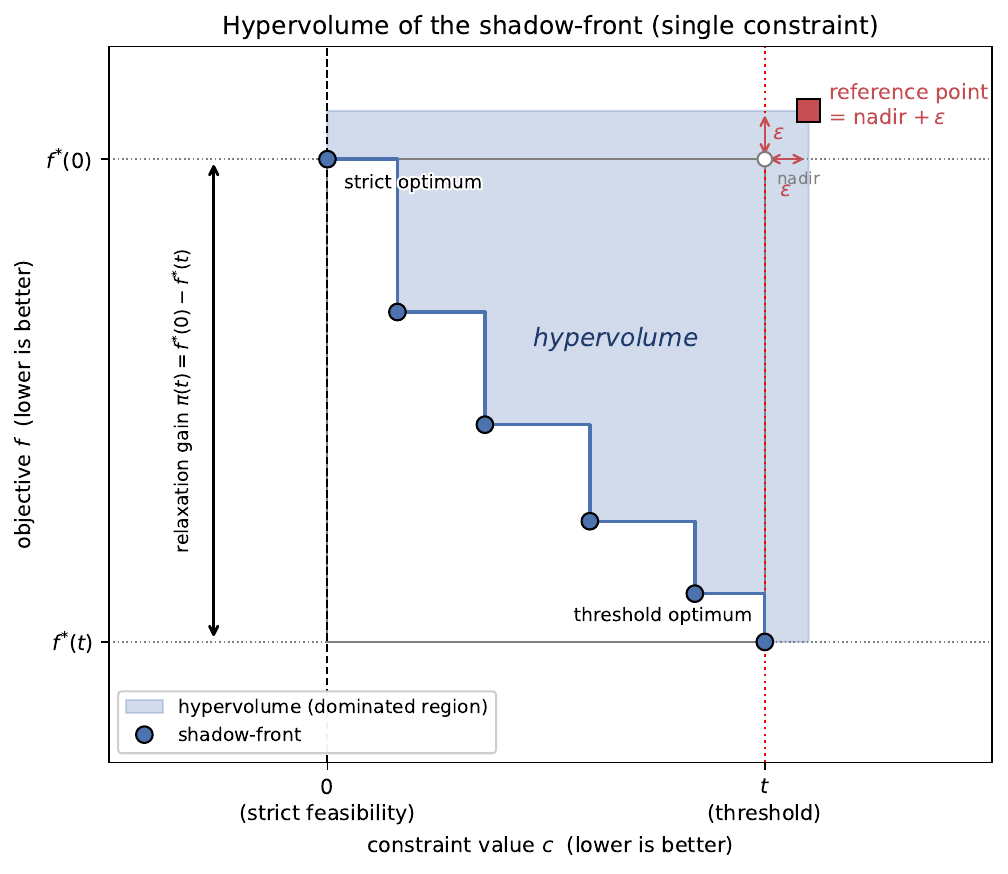}
      \caption{Shadow-front and its hypervolume for a single constraint. As the
            constraint is relaxed from its original boundary ($c=0$) up to the threshold
            ($c=t$), the best attainable objective improves from $f^\star(0)$ to
            $f^\star(t)$; the non-dominated trade-off between relaxation and objective is
            the shadow-front, and the vertical extent $f^\star(0) - f^\star(t)$ is the
            relaxation gain $\pi(t)$.
            The hypervolume of the front with respect to the nadir, expanded by a small
            $\varepsilon$ on each axis (exaggerated here), measures how well the front has
            been mapped out.}
      \label{fig:hv_illustration}
\end{figure}

Taking the thresholds as bounds rather than fixed relaxations leads to the
\emph{shadow-front}:\ the set of non-dominated points
$(r_1(\mathbf{x}),\allowbreak \dots,\allowbreak r_k(\mathbf{x}),\allowbreak f(\mathbf{x}))$ over
threshold-feasible $\mathbf{x}$, where $r_j(\mathbf{x}) = \max(0,
      c_j(\mathbf{x}))$ is the relaxation that $\mathbf{x}$ uses on constraint $j$,
and all coordinates are minimized. The shadow-front is the Pareto front of
relaxation against objective value restricted to the box $[0, t_1] \times \dots
      \times [0, t_k]$. Its two end points are the strict optimum, with $r_j = 0$ for
all $j$, and the threshold optimum $f^\star(\mathbf{t})$, and their objective
difference is the relaxation gain $\pi(\mathbf{t})$ of the full relaxation.
Figure~\ref{fig:hv_illustration} visualizes the single-constraint case.

What a decision maker may ask from these objects varies along two dimensions,
which give the four scenarios in Table~\ref{tab:scenarios}. The first dimension
is \emph{end points versus curve}. A decision maker who has already settled on
the relaxation they can afford needs only the gain $\pi(\mathbf{t})$ at the end
point; one who is still deciding how far to go needs the whole shadow-front, to
see where additional relaxation stops paying off. The second dimension is
\emph{single versus joint relaxation}. Relaxing one constraint at a time, that
is, $\mathbf{s} = s_j \mathbf{e}_j$, where $\mathbf{e}_j$ is the $j$-th unit
vector of $\mathbb{R}^k$, with the other constraints held at their
original boundary, attributes the gain to that constraint alone, which is the
spirit of a shadow price and the natural question when each constraint has its
own owner or cost. Relaxing all constraints jointly, $\mathbf{s} \leq
      \mathbf{t}$ component-wise, answers what the best achievable outcome is when
every limit may move, as in a negotiation. Under single relaxation the
shadow-front decomposes into $k$ one-dimensional curves, one per constraint,
each with its own gain $\pi(t_j \mathbf{e}_j)$; under joint relaxation it is a
single $k$-dimensional surface with the joint gain $\pi(\mathbf{t})$ at its far
end. These scenarios are not exclusive; a decision maker may well ask for the
per-constraint gains first and the joint surface afterwards, or for the end point of one constraint and the shadow-front of another.

\begin{table}[htbp]
      \centering
      \caption{The four decision-maker scenarios and the selection modes that
	address them. The relaxed mode targets the end points of each
	single-constraint curve, the joint mode targets the end point when all
	constraints are relaxed simultaneously, the ranking mode spreads over
	the single-constraint curves of the relaxed mode, and the hypervolume
	mode spreads over the joint surface.}
      \label{tab:scenarios}
      \begin{tabular}{@{}lll@{}}
            \toprule
                        & \multicolumn{2}{l}{Relaxation}                    \\
            \cmidrule(l){2-3}
            Target      & Single                         & Joint            \\
            \midrule
            End points  & Relaxed mode                   & Joint mode \\
            Whole curve & Ranking mode                   & Hypervolume mode \\
            \bottomrule
      \end{tabular}
\end{table}

Each scenario calls for a different emphasis in the search:\ on the end points
or on the spread of the front, and on single relaxation directions or on the
joint relaxation box. The four selection modes we propose below are designed
along these lines, and Table~\ref{tab:scenarios} shows which mode addresses
which scenario. Each cell has a dedicated mode; the hypervolume mode, which
protects the two end points of the joint shadow-front, also serves the end
points under joint relaxation. The thresholds are meant to come from someone familiar with the problem, such
as the decision maker, who states the relaxation they are comfortable with. In
this study we set them ourselves (Section~\ref{subsec:exp_setup}). Our claim is
that, with a suitably designed environmental selection operator, the
information required in each scenario is obtained as a by-product of the
evolutionary algorithm that finds the optimum, at no cost in additional
objective evaluations.

\subsection{Selection operator and baseline mode}
\label{subsec:selection_op}

\begin{figure}[!htb]
      \centering
      \includegraphics[width=\linewidth]{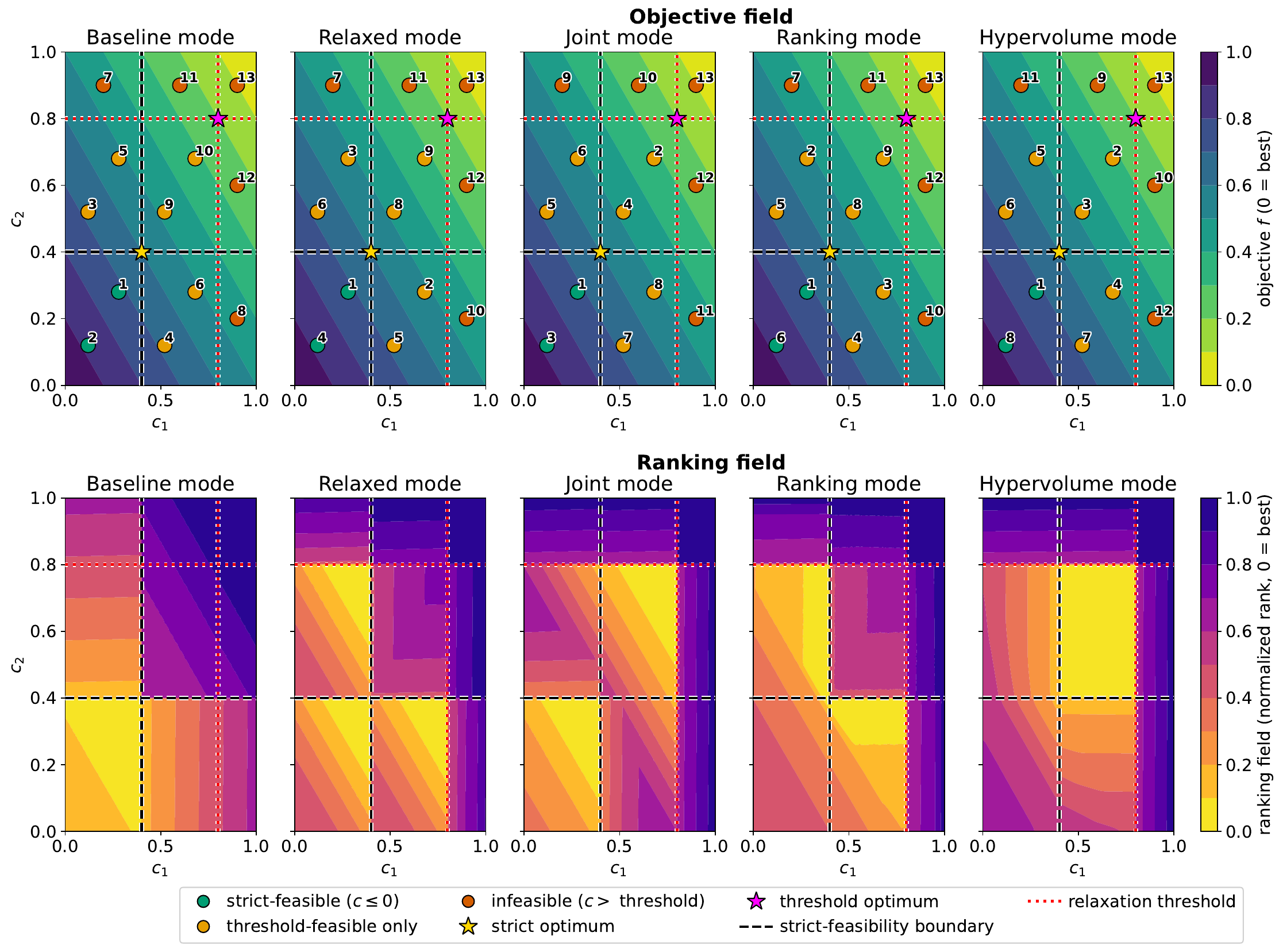}
      \caption{How each mode ranks a common set of solutions on a synthetic
            two-constraint landscape. Both constraints and the objective are scaled to
            $[0,1]$; the objective improves along the diagonal from the bottom-left to the
            top-right corner. The strict-feasibility boundary (dashed) and the relaxation
            threshold (dotted) split each axis into three bands. Top row:\ the solutions,
            colored by feasibility category, annotated with the rank each mode assigns
            ($1 =$ best). Bottom row:\ the normalized rank of each point when a dense grid
            of solutions covering the square is ranked together ($0 =$ best). The baseline ranks every strictly feasible solution
            above every other solution, whatever its objective. The relaxed, joint, and
            ranking modes also rank well solutions in the relaxed bands, between the
            strict-feasibility boundary and the threshold. The hypervolume mode ranks
            highest the solutions that relax both constraints, in the box between the
            strict-feasibility boundary and the threshold; inside the box its ranks
            alternate between neighboring points, which we smooth within each region for
            display.}
      \label{fig:mode_illustration}
\end{figure}

We consider a generational evolutionary algorithm in which parents and
offspring are pooled and the pool is ranked and truncated to the population
size. Our selection operator comes in a baseline mode and four alternative
modes. Figure~\ref{fig:mode_illustration} illustrates how each ranks a common
set of solutions on a synthetic two-constraint landscape.

All modes share one ranking rule. Given a violation vector
$v_j(\mathbf{x})$, $j=1,\dots,k$, solutions with all $v_j = 0$ are ranked first,
in ascending order of $f$. The remaining solutions follow, sorted
lexicographically by the number of violated constraints, then by the sum of
normalized violations $\sum_j \bar{v}_j(\mathbf{x})$, and then by $f$. Here
$\bar{v}_j$ is $v_j$ min-max scaled to $[0,1]$ over the infeasible solutions in
the pool (and set to $0$ when the range is zero), so that constraints of
different magnitudes contribute comparably. Sorting first by the number of
violated constraints prefers solutions that are closer to feasibility in a
scale-free sense. We call this
the \emph{feasibility-first lexicographic rule}. It is applied with different
violation vectors by the different modes:\ with the original violations
$v_j(\mathbf{x}) = \max(0, c_j(\mathbf{x}))$ it combines the feasibility rule
of \citet{deb2000efficient}, which prefers feasible solutions and orders them by
objective, with the comparison of infeasible solutions
of \citet{coello2000constrainthandling}, first by the number of violated
constraints and only in case of a tie by the amount of violation. Unlike
in~\citet{deb2000efficient}, where the constraints are normalized in advance,
we scale the violations over the infeasible part of the pool, and we break
the remaining ties by $f$. With constraint $j$ relaxed it uses
$v_j(\mathbf{x}) = \max(0, c_j(\mathbf{x}) - t_j)$. Before any ranking, duplicate
solutions, that is, solutions with identical decision variable vectors, are
removed from the pool to maintain diversity.

The \emph{baseline mode} ranks the pool once, using the original violations, and
truncates it to the population size (Algorithm~\ref{alg:baseline_mode}). It
represents conventional constraint handling that produces no relaxation
information and serves as the control in our experiments.

\begin{algorithm}[t]
      \caption{Baseline Mode}
      \label{alg:baseline_mode}

      \KwIn{parents $P$, offspring $Q$, objective $f$, constraints $c_1,\dots,c_k$, population size $N$}
      \KwOut{new population $P'$}

      $R \leftarrow P \cup Q$\;
      remove duplicate solutions from $R$\;

      partition $R$ into feasible ($c_j \le 0\ \forall j$) and infeasible\;

      rank feasible solutions by $f$ (ascending)\;

      rank infeasible solutions lexicographically by\\
      \hspace{1em}(number of violated constraints,\ normalized violation sum,\ $f$)\;

      $O \leftarrow$ ranked feasible followed by ranked infeasible\;

      $P' \leftarrow$ first $N$ solutions from $O$\;

      \Return{$P'$}\;
\end{algorithm}

\FloatBarrier
\subsection{Relaxed mode:\ end points under single relaxation}
\label{subsec:relaxed_mode}

The \emph{relaxed mode} ranks the pool $k+1$ times. Ranking $O_0$ uses the
original constraints. Ranking $O_j$, for $j = 1,\dots,k$, relaxes constraint $j$
alone to $c_j \leq t_j$ and holds every other constraint at its original
boundary. The new population is then filled by round-robin selection over
$O_0, \dots, O_k$:\ the highest-ranked remaining solution is taken from each
ranking in turn, and a solution taken from one ranking is removed from all of
them (Algorithm~\ref{alg:relaxed_mode}).

The mode is built for the end points under single relaxation. The top of $O_0$
pursues the strict optimum, and the top of $O_j$ pursues the optimum with
constraint $j$ relaxed, so the population always carries the best known
estimate of each per-constraint gain $\pi(t_j \mathbf{e}_j)$. Relaxing a single
constraint per ranking keeps the construction faithful to the shadow price:\
the improvement shown by the top solutions of $O_j$ is attributable to
constraint $j$ alone. It also keeps each ranking focused on the boundary region
of its own constraint, since the top of $O_j$ consists of solutions that
satisfy all other constraints and lie in the relaxed band $0 < c_j \leq t_j$,
which are exactly the solutions that reveal what relaxing $c_j$ buys. The
construction scales linearly in $k$. It exerts no
pressure to spread solutions along the curve between the end points.

\begin{algorithm}[t]
      \caption{Relaxed Mode}
      \label{alg:relaxed_mode}

      \KwIn{parents $P$, offspring $Q$, objective $f$, constraints $c_1,\dots,c_k$, thresholds $t_1,\dots,t_k$, population size $N$}
      \KwOut{new population $P'$}

      $R \leftarrow P \cup Q$\;

      remove duplicate solutions from $R$\;

      $O_0 \leftarrow$ rank $R$ using strict feasibility $c_j \le 0$ for $j=1,\dots,k$\;

      \For{$j = 1,\dots,k$}{
            $O_j \leftarrow$ rank solutions in $R$ using relaxed feasibility $c_j \le t_j$ and $c_i \le 0$ for all $i \neq j$\;
      }

      $P' \leftarrow \emptyset$\;

      \While{$|P'| < N$}{
            select next non-empty pool index in round-robin order $j\in\{0,\dots,k\}$\;
            select the highest-ranked remaining solution from $O_j$\;
            remove it from all pools $O_0,\dots,O_k$\;
            add it to $P'$\;
      }

      \Return{$P'$}\;

\end{algorithm}

\FloatBarrier
\subsection{Joint mode:\ end points under joint relaxation}
\label{subsec:joint_mode}

The \emph{joint mode} is the relaxed mode with two rankings instead of $k+1$.
Ranking $O_0$ uses the original constraints, as before, and ranking
$O_{\mathbf{t}}$ relaxes every constraint to its threshold at the same time,
$c_j \leq t_j$ for all $j$. The new population is filled by round-robin
selection over $O_0$ and $O_{\mathbf{t}}$, as in the relaxed mode
(Algorithm~\ref{alg:joint_mode}).

The mode is built for the end points under joint relaxation. The top of $O_0$
pursues the strict optimum and the top of $O_{\mathbf{t}}$ the threshold
optimum, so the population always carries the best known estimates of both
ends of the joint shadow-front, whose difference is the joint gain
$\pi(\mathbf{t})$. In exchange, it gives up the attribution of the gain to
individual constraints, which the per-constraint rankings of the relaxed mode
keep, and, like the relaxed mode, it exerts no pressure to spread solutions
between the end points. With a single constraint its two rankings are those of
the relaxed mode, and the two modes coincide. Since it ranks the pool only
twice, its cost does not grow with the number of constraints beyond computing
the violations.

\begin{algorithm}[t]
      \caption{Joint Mode}
      \label{alg:joint_mode}

      \KwIn{parents $P$, offspring $Q$, objective $f$, constraints $c_1,\dots,c_k$, thresholds $t_1,\dots,t_k$, population size $N$}
      \KwOut{new population $P'$}

      $R \leftarrow P \cup Q$\;

      remove duplicate solutions from $R$\;

      $O_0 \leftarrow$ rank $R$ using strict feasibility $c_j \le 0$ for $j=1,\dots,k$\;

      $O_{\mathbf{t}} \leftarrow$ rank $R$ using relaxed feasibility $c_j \le t_j$ for $j=1,\dots,k$\;

      $P' \leftarrow \emptyset$\;

      \While{$|P'| < N$}{
            select the next non-empty pool in round-robin order $O_0, O_{\mathbf{t}}$\;
            select the highest-ranked remaining solution from it\;
            remove it from both pools\;
            add it to $P'$\;
      }

      \Return{$P'$}\;

\end{algorithm}

\FloatBarrier
\subsection{Ranking mode:\ whole curve under single relaxation}
\label{subsec:ranking_mode}

The \emph{ranking mode} forms the same $k+1$ rankings as the relaxed mode but
combines them differently. Each solution is assigned a $(k+1)$-dimensional rank
vector holding its position in $O_0, \dots, O_k$. These rank vectors are treated
as the objectives of a multiobjective problem and subjected to non-dominated
sorting; within each front, the mean distance to the $\kappa$ nearest neighbors in
rank-vector space serves as the niching measure. A solution's fitness is its front index
plus a niching term in $[0,1)$, and the $N$ fittest solutions survive
(Algorithm~\ref{alg:ranking_mode}).

The mode is built for the whole curve under single relaxation. Where the
relaxed mode interleaves independent picks, the ranking mode retains solutions
that rank well jointly across the original and all single-constraint relaxed
orderings, and the niching spreads them along the trade-off between the
strict ranking and each relaxed ranking, that is, along each per-constraint
curve of the shadow-front. The hypothesis behind this design is that
recovering this surrounding trade-off structure, rather than a sequence of
single relaxed points, also yields a more reliable estimate of the end points.

\begin{algorithm}[t]
      \caption{Ranking Mode}
      \label{alg:ranking_mode}

      \KwIn{parents $P$, offspring $Q$, objective $f$, constraints $c_1,\dots,c_k$, thresholds $t_1,\dots,t_k$, population size $N$}
      \KwOut{new population $P'$}

      $R \leftarrow P \cup Q$\;

      remove duplicate solutions from $R$\;

      construct ranking pools $O_0,\dots,O_k$ as in Relaxed Mode\;

      \ForEach{$x \in R$}{
            form rank vector\\
            \hspace{1em}$r(x) = ( \text{position in } O_0,\dots,\text{position in } O_k )$
      }

      apply non-dominated sorting to $\{r(x)\}$\;

      compute the mean distance to the $\kappa$ nearest neighbors within each front\;

      assign fitness based on front rank and this mean distance\;

      $P' \leftarrow$ $N$ solutions with best fitness\;

      \Return{$P'$}\;
\end{algorithm}

\FloatBarrier
\subsection{Hypervolume mode:\ joint relaxation}
\label{subsec:hypervolume_mode}

The relaxed and ranking modes relax one constraint at a time, and the joint
mode relaxes all constraints together but exerts no pressure to spread
solutions over the joint relaxation box. The
\emph{hypervolume mode} does both. It treats all relaxations jointly and ranks
threshold-feasible solutions by their contribution to the hypervolume of the
shadow-front.

Each threshold-feasible solution is mapped to the point
$(f(\mathbf{x}),\allowbreak r_1(\mathbf{x}),\allowbreak \dots,\allowbreak r_k(\mathbf{x}))$ of
Section~\ref{subsec:definitions}, with every axis normalized to $[0,1]$:\ the
objective by its range over the current threshold-feasible pool, and each
relaxation $r_j$ by its threshold $t_j$. The reference point is the upper
corner of the normalized box (Figure~\ref{fig:hv_illustration}), expanded outward by a small fraction
$\varepsilon$ on every axis, that is, $1 + \varepsilon$ on every axis, so that
solutions on the boundary of the box still contribute volume; the
normalization is recomputed from each pool. As in SMS-EMOA~\citep{beume2007sms}, the points
are first sorted into non-dominated fronts. Within each front the solutions are then ordered
greedily:\ the solution whose removal reduces the hypervolume of the front
least is removed, contributions are recomputed, and the process is repeated,
so that the survival order reflects each solution's marginal value to the
shrinking set. Fronts are ordered only until $N$ solutions are placed. The two end
points of the front, the best strictly feasible and the best threshold-feasible
solution, both of which lie on the first front, are protected from removal
until last; their difference is the joint relaxation gain $\pi(\mathbf{t})$
itself. Solutions
beyond the thresholds are placed after the threshold-feasible ones, ordered by
the feasibility-first lexicographic rule with violations measured from the
thresholds (Algorithm~\ref{alg:hypervolume_mode}).

The hypervolume mode thus serves both joint-relaxation scenarios of
Table~\ref{tab:scenarios}:\ the whole surface through contribution-based
selection, and its end points through the protection. Any mode that finds the
joint end point must evaluate solutions inside the relaxation box on the way;
spreading those solutions over the box, as the hypervolume mode does, delivers
the surface for the same evaluations. The joint mode
(Section~\ref{subsec:joint_mode}) targets the end points alone, at a much lower
selection cost.

\begin{algorithm}[t]
      \caption{Hypervolume Mode}
      \label{alg:hypervolume_mode}

      \KwIn{parents $P$, offspring $Q$, objective $f$, constraints $c_1,\dots,c_k$, thresholds $t_1,\dots,t_k$, population size $N$, reference margin $\varepsilon$}
      \KwOut{new population $P'$}

      $R \leftarrow P \cup Q$\;

      remove duplicate solutions from $R$\;

      partition $R$ into threshold-feasible $T$ ($c_j \le t_j$ for all $j$) and beyond-threshold $B$\;

      \ForEach{$x \in T$}{
            form the point $y(x) = ( f(x),\ r_1(x),\dots,r_k(x) )$,\ where $r_j(x) = \max(0, c_j(x))$\;
      }

      normalize each axis of $\{y(x)\}$ to $[0,1]$:\ the objective by its range over $T$, each $r_j$ by $t_j$\;

      $z \leftarrow (1+\varepsilon, \dots, 1+\varepsilon)$, the upper corner of the normalized box expanded by $\varepsilon$\;

      sort $T$ into non-dominated fronts $F_1, F_2, \dots$ of the normalized points\;

      mark the best-objective and the best strictly-feasible solution of $F_1$ as protected\;

      $O_T \leftarrow$ empty list\;

      \ForEach{front $F_i$, in order, while $|O_T| < N$}{
            $S \leftarrow F_i$;\quad $L \leftarrow$ empty list (removal order, worst first)\;
            \While{$|S| > 1$}{
                  $x^\star \leftarrow$ least-contributing solution in $S$ to its hypervolume w.r.t.\ $z$,\\
                  \hspace{1em}skipping protected solutions unless only they remain\;
                  append $x^\star$ to $L$;\ remove $x^\star$ from $S$\;
            }
            append the remaining solution of $S$ to $L$\;
            append the reverse of $L$ (best to worst) to $O_T$\;
      }

      $O_B \leftarrow$ rank $B$ by the feasibility-first lexicographic rule\;

      $P' \leftarrow$ first $N$ solutions of $O_T$ followed by $O_B$\;

      \Return{$P'$}\;
\end{algorithm}

\FloatBarrier
\section{Experiments and results}
\label{sec:experiments}

\subsection{Setup}
\label{subsec:exp_setup}

We embed the selection operator in the generational scheme of
NSGA-II~\citep{deb2002fast}, with the variation operators and parameter values
of DESDEO's standard NSGA-II settings.
The initial population is drawn by Latin hypercube sampling. In every
generation, the selection operator ranks the pool formed by the current
population and its offspring and keeps the best $N$ solutions as the new
population. From this population, $N$ parents are then chosen by binary
tournaments drawn with replacement. The parents are shuffled and paired,
and each pair is recombined with probability $0.9$ by simulated binary
crossover~\citep{deb1995simulated}. Within each pair, variables are
crossed with probability $0.5$, using the bounded form of the operator, whose
spread distribution (distribution index $\eta_c = 20$) is truncated so that the
children stay within the variable bounds~\citep{deb2001truss}, and the two
children's values of a crossed variable are exchanged with probability $0.5$.
The other variables are inherited unchanged. Each variable of a child is then
mutated with probability $1/n$ by bounded polynomial
mutation~\citep{deb2001multiobjective} (distribution index $\eta_m = 20$), and
the children are clipped to the bounds; together with the population they form
the next pool. The
ranking mode uses as its niching measure the mean distance of a solution to its
three nearest neighbors ($\kappa=3$), and the hypervolume mode a reference
margin of $\varepsilon = 0.01$ (Figure~\ref{fig:hv_illustration}). Each run
lasts $200$ generations. We consider small population sizes, $N \in \{12,
      36, 60\}$, to mimic the evaluation budgets of expensive black-box problems,
since our claim is that the shadow-front is obtained at such budgets without
extra evaluations. The three population sizes thus amount to three evaluation
budgets per run, $(200+1) \times N$ ($N$ evaluations for the initial population and $N$ in each
of the $200$ generations), that is, $2\,412$,
$7\,236$, and $12\,060$  evaluations.

We use ten benchmark problems with one to four constraints and two to twenty
variables (Table~\ref{tab:problems_thresholds}), all accessed through the DESDEO
framework~\citep{misitano2021desdeo}. Branin and Mystery
follow~\citet{sasena2002flexibility} and Townsend follows~\citet{townsend2014}.
The cantilevered beam~\citep{hexaly2025cantilever}, the pressure vessel~\citep{kannan1994augmented}, and
G02, G06, G08, G09, and G24 from the CEC~2006 suite~\citep{liang2006problem} are
provided through DESDEO's interface to pymoo~\citep{blank2020pymoo} and are
used as implemented in pymoo (version 0.6.2). The source code used to run the experiments is archived on Zenodo~\citep{misitano2026code}.

For each problem we need reference values of the relaxed optima against which
to judge the modes:\ the threshold optimum $f^\star(\mathbf{t})$ and the
single-relaxed optima $f^\star(t_j \mathbf{e}_j)$ at every level. The problems
are nonconvex, so we use best-known values combined from three sources. First,
a reference front obtained by treating the $k$ constraints as additional
objectives and running NSGA-III~\citep{deb2014nsga3} with a large budget
($10\,000$ generations, population size $200$). We use DESDEO's standard
NSGA-III settings:\ $200$ reference directions on the simplex, the variation
operators described above with every pair recombined and $\eta_c = 30$, and an
initial population that consists of the known optimum, perturbed copies of it
for $20$ percent of the population (Gaussian, with a standard deviation of $2$
percent of each variable's range), and uniformly random points. With $2 \times 10^6$
evaluations, this budget is more than $150$ times the largest budget of the
experiment runs, which makes the front a strong independent reference. Second, the best
threshold-feasible solutions of the experiment runs themselves, pooled over all
modes, threshold levels, and population sizes. Third,
multistart sequential quadratic programming (SLSQP) on each relaxed problem,
started from the best run solutions and from random points in the box. The
best feasible value of the three is the reference, and together with the known
optimum $f^\star(\mathbf{0})$ it gives the reference gains
via~\eqref{eq:shadow_price}. The front alone is not enough:\ on G06, G09, and
the pressure vessel it falls short of the reference gain by up to $58$ percent
at the joint end points, the local solver mostly improves markedly on it, and most runs
beat it on G09 and on the pressure vessel at the low and medium levels; on G02
the local solver finds a small gain that the front misses entirely. Since the
references include the runs, no run can be better than them:\ a gap of zero
means that a run has found the best-known value. The same procedure applied
to the strict problem recovers the known optima, with one exception:\ on the
pressure vessel it finds $5885.33$, because the commonly cited $6059.72$ holds
for shell thicknesses restricted to multiples of $0.0625$, whereas the pymoo
implementation treats them as continuous. We use the continuous optimum, which
makes the results comparable to other studies that use the pymoo
implementation of the problem.

Each problem has three threshold levels (low, medium, high;
Table~\ref{tab:problems_thresholds}). We first derived candidate values from the
reference front, from the violations of front solutions that improve on the
known optimum. Several of these were, however, unattainable when a constraint
is relaxed on its own:\ a front solution may improve on the optimum only by
relaxing several constraints jointly, so the violation it shows on one
constraint cannot be reached by relaxing that constraint alone.
We therefore inspected each constraint over the design space with contour plots
and adjusted the thresholds so that every level is a meaningful, attainable
relaxation. A constraint that is not binding and cannot be meaningfully relaxed is held
at its original boundary. Even so, relaxing a constraint on its own does not
always improve the optimum, because another constraint binds first. With the
best-known references this is the case, at every level, for $c_2$ of the
cantilevered beam, $c_2$ of G02, $c_1$ of G08, and $c_2$ and $c_3$ of G09.
These constraints still take part in the joint relaxation, but the
single-relaxation metrics of Section~\ref{subsec:metrics} leave them out. A
decision maker cannot know in advance which constraints these are, and does not
need to:\ in our runs, such a constraint shows no gain at all, which is the
correct answer, except on G08, where $7$ of $4\,500$ runs report a small gain
for $c_1$ of up to $5.5$ percent of the optimum, an artifact of not finding the correct constrained optimum. The metrics leave them out only because their
relative gap, normalized by a reference gain of zero, is undefined. On
G02 the attainable gain is tiny:\ relaxing $c_1$ improves the optimum by
$0.3$, $0.6$, and $1.0$ percent of its magnitude at the three levels. We keep
these problems and constraints deliberately instead of selecting benchmarks
on which every constraint pays off on its own:\ constraints that do not bind,
that bind only together with others, or whose relaxation gains little are
features that real problems are likely to exhibit as well, and a method should
not report gains where there are none.

\begin{table}[htbp]
      \centering
      \caption{Benchmark problems and constraint relaxation thresholds. For each
            problem:\ number of variables, number of constraints, known optimum, and the low/medium/high
            threshold for every constraint. $^{\dagger}$ marks a non-binding constraint held
            at its original boundary and not relaxed.}
      \label{tab:problems_thresholds}
      \small
      \setlength{\tabcolsep}{4pt}
      \begin{tabular}{l c c r l r r r}
            \toprule
            \textbf{Problem}               & \textbf{\#Var}      & \textbf{\#Con}     & \textbf{Optimum}               & \textbf{Constr.} & \textbf{Low}                               & \textbf{Med.} & \textbf{High} \\
            \midrule
            Branin                         & 2                   & 1                  & $-268.78792$                   & $c_1$            & 15.0                                       & 40.0          & 60.0          \\
            \midrule
            Mystery                        & 2                   & 1                  & $-1.174261$                    & $c_1$            & 0.15                                       & 0.25          & 0.35          \\
            \midrule
            \multirow{2}{*}{Cant.\ beam}   & \multirow{2}{*}{4}  & \multirow{2}{*}{2} & \multirow{2}{*}{$92.7693$}
                                           & $c_1$               & 0.5                & 1.0                            & 1.5                                                                                           \\
                                           &                     &                    &                                & $c_2$            & 0.65                                       & 0.95          & 1.4           \\
            \midrule
            \multirow{4}{*}{Pres.\ vessel} & \multirow{4}{*}{4}  & \multirow{4}{*}{4} & \multirow{4}{*}{$5885.33275$}
                                           & $c_1$               & 0.06329            & 0.1266                         & 0.2532                                                                                        \\
                                           &                     &                    &                                & $c_2$            & 0.03076                                    & 0.06152       & 0.1230        \\
                                           &                     &                    &                                & $c_3$            & 0.04972                                    & 0.09943       & 0.1989        \\
                                           &                     &                    &                                & $c_4^{\dagger}$  & \multicolumn{3}{c}{inactive (not relaxed)}                                 \\
            \midrule
            Townsend                       & 2                   & 1                  & $-2.0239884$                   & $c_1$            & 0.07450                                    & 0.1490        & 0.2980        \\
            \midrule
            \multirow{2}{*}{G02}           & \multirow{2}{*}{20} & \multirow{2}{*}{2} & \multirow{2}{*}{$-0.8036191$}
                                           & $c_1$               & 0.05               & 0.10                           & 0.15                                                                                          \\
                                           &                     &                    &                                & $c_2$            & 20.0                                       & 40.0          & 60.0          \\
            \midrule
            \multirow{2}{*}{G06}           & \multirow{2}{*}{2}  & \multirow{2}{*}{2} & \multirow{2}{*}{$-6961.81388$}
                                           & $c_1$               & 0.25               & 0.50                           & 0.75                                                                                          \\
                                           &                     &                    &                                & $c_2$            & 0.15                                       & 0.25          & 0.35          \\
            \midrule
            \multirow{2}{*}{G08}           & \multirow{2}{*}{2}  & \multirow{2}{*}{2} & \multirow{2}{*}{$-0.09582504$}
                                           & $c_1$               & 0.25               & 0.50                           & 0.75                                                                                          \\
                                           &                     &                    &                                & $c_2$            & 10.0                                       & 15.0          & 20.0          \\
            \midrule
            \multirow{4}{*}{G09}           & \multirow{4}{*}{7}  & \multirow{4}{*}{4} & \multirow{4}{*}{$680.630057$}
                                           & $c_1$               & 100.0              & 200.0                          & 300.0                                                                                         \\
                                           &                     &                    &                                & $c_2$            & 2.0                                        & 4.0           & 6.0           \\
                                           &                     &                    &                                & $c_3$            & 2.0                                        & 4.0           & 6.0           \\
                                           &                     &                    &                                & $c_4$            & 10.0                                       & 15.0          & 20.0          \\
            \midrule
            \multirow{2}{*}{G24}           & \multirow{2}{*}{2}  & \multirow{2}{*}{2} & \multirow{2}{*}{$-5.50801327$}
                                           & $c_1$               & 0.09880            & 0.1976                         & 0.3952                                                                                        \\
                                           &                     &                    &                                & $c_2$            & 0.2                                        & 0.4           & 0.8           \\
            \bottomrule
      \end{tabular}
\end{table}

In total we run the baseline and the four alternative modes, three threshold
levels, and three population sizes on ten problems, totaling $450$
configurations, each with $100$ independent
runs. Each run records every solution it evaluates, with its objective and
constraint values and the generation in which it was evaluated; we call this
record the run's \emph{archive}. All metrics of Section~\ref{subsec:metrics} are
computed from it, at each generation from the solutions evaluated up to that
generation. Every run is seeded, and
the $r$-th run of every configuration uses the same seed, so that all five
modes start each run from the exact same initial population and are compared on
common random numbers. Since paired tests on $100$ runs detect small
differences, we report on significance as well as the size of the
differences. The whole pipeline, from problem definitions to the tables and
figures in this paper, is implemented in DESDEO~\citep{misitano2021desdeo}
(version 2.9.0) and orchestrated with Snakemake~\citep{molder2021snakemake}. The
source code, exact software versions, and instructions~\citep{misitano2026code}
and the run data, reference fronts, and derived results~\citep{misitano2026data}
are archived on Zenodo. Each data file records the complete algorithm settings it was produced
with, and rerunning a configuration reproduces its data exactly, so that the
experiments can be checked or extended to other problems.

\subsection{Performance metrics}
\label{subsec:metrics}

We report five metrics per generation of each run:\ the optimization
performance, and one metric for each cell of Table~\ref{tab:scenarios}, so
that every mode is judged on the scenario it was designed for as well as on
the others. The reference values are the known optimum $f^\star(\mathbf{0})$
and the best-known relaxed optima of Section~\ref{subsec:exp_setup}, and every
scenario metric is normalized so that
its values are comparable across problems and numbers of constraints.

\paragraph{Best-so-far feasible objective.} The cumulative minimum over
generations of the best strictly feasible objective value ($c_j \leq 0$ for all
$j$). This measures standard optimization performance.

\paragraph{End points under single relaxation:\ relative gap.} For constraint
$j$, consider the solutions found so far that satisfy $c_j \leq t_j$ and every
other constraint strictly. Their best objective value $\hat f_j$ estimates the
single-relaxed optimum $f^\star(t_j \mathbf{e}_j)$, and we report its relative
optimality gap
\begin{equation}
      \label{eq:single_gap}
      \gamma_j = \frac{\max\{0,\; \hat f_j - f^\star(t_j \mathbf{e}_j)\}}{\pi(t_j \mathbf{e}_j)}.
\end{equation}
The gap is $0$ when the run has found the relaxed optimum and $1$ when it has
found nothing better than the strict optimum. It is clipped at zero only to guard
against rounding, since the reference includes the runs and no run improves on it. The metric is the mean of $\gamma_j$ over the constraints that have a
positive reference gain under single relaxation. A constraint whose relaxation
alone does not improve the optimum, because another constraint binds first,
carries no information about the end points and is left out
(Section~\ref{subsec:exp_setup}). Since the gain is attributed to one
constraint at a time, this metric measures how well a mode recovers the shadow
price of each constraint at its threshold.

\paragraph{Whole curve under single relaxation:\ hypervolume.} For the same
solutions we take, per run, the non-dominated points of the archive
$(r_j(\mathbf{x}), f(\mathbf{x}))$ and compute its hypervolume in the box
$[f^\star(t_j \mathbf{e}_j), f^\star(\mathbf{0})] \times [0, t_j]$. Each axis
is normalized to $[0,1]$ and the reference point is the upper corner of the
box, the nadir, expanded by $\varepsilon = 0.01$ on each normalized axis, that is,
by $1\%$ of the axis range, so that solutions
on the boundary of the box, in particular the two end points, still contribute
volume (Figure~\ref{fig:hv_illustration}). Before computing the hypervolume,
solutions are mapped into the box:\ a strictly satisfied constraint coordinate
is clipped to $0$, an objective value better than the reference to the
objective floor, and solutions outside the box are discarded. The metric is
the mean over the same contributing constraints. It rewards both how far a
mode pushes the attainable gain of each constraint and how evenly it covers the
relaxation range in between.

\paragraph{End points under joint relaxation:\ relative gap.} The best-so-far
threshold-feasible objective value ($c_j \leq t_j$ for all $j$) estimates the
threshold optimum $f^\star(\mathbf{t})$, and we report its relative optimality
gap, defined as in~\eqref{eq:single_gap} with $f^\star(\mathbf{t})$ and
$\pi(\mathbf{t})$ in place of the single-relaxation quantities. This measures
how well a mode recovers the far end point of the joint shadow-front.

\paragraph{Whole curve under joint relaxation:\ hypervolume.} The hypervolume
of the non-dominated threshold-feasible solutions of the run's archive in the box
$[f^\star(\mathbf{t}), f^\star(\mathbf{0})] \times [0, t_1] \times \dots
      \times [0, t_k]$, with every axis normalized to $[0,1]$, strictly satisfied
constraint coordinates clipped to $0$, objective values better than the reference
clipped to the objective floor, and the reference point at $1+\varepsilon$ on every axis. Constraints that
are not relaxed (zero threshold) are excluded from the box but must still be
satisfied. Only threshold-feasible solutions count; the margin $\varepsilon$
expands the reference point only, not the thresholds.

Note that when measuring performance by hypervolume the hypervolume mode has an advantage by construction, as it uses the performance metric inside. For this reason, in the results, we additionally compare the remaining modes among themselves, as if the hypervolume mode had not been run:\ the same rule as for the winners, applied to the other four modes with the Holm correction over their $6$ pairwise tests.

On G02 the reference gain is below one percent of the optimum, but within the
budget the runs end far from the optimum of this $20$-dimensional problem:\ on
average, their best strictly feasible values are $15$ to $52$ percent worse than
the strict optimum.  Its
hypervolumes are zero in all but one configuration, in which a single run
reaches part of the joint gain. We keep it in the tables as a sanity check:\
no mode reports a gain that is not there.

For each configuration we aggregate every metric over the $100$ runs by its
mean and a $95\%$ Student's $t$ confidence interval per generation. Generations
in which a run has not yet found a solution of the required kind contribute no
value; this affects mainly early generations, although on G06 at $N = 12$ one
or two runs of some modes find no strictly feasible solution within the budget.
For every problem and metric we
plot the mean and confidence interval against the number of evaluations, for
each combination of population size and threshold level, in the supplementary
material.

\FloatBarrier
\subsection{Results}
\label{subsec:observations}

Tables~\ref{tab:best-so-far} to~\ref{tab:hv} summarize the final-generation
means of the five metrics over the $100$ runs of all $450$ configurations:\
one table for the best-so-far objective and one for each cell of
Table~\ref{tab:scenarios}. Each cell holds the mean for the baseline (Bs),
relaxed (Rx), joint (Jt), ranking (Rk), and hypervolume (Hv) modes. In problems
with only a single constraint, the relaxed and joint modes will behave identically.
Significance is
assessed with the 10 pairwise two-sided Wilcoxon signed-rank tests of a cell
at $\alpha = 0.05$, with the Holm correction over the 10 tests~\citep{holm1979simple}. Bold marks the winning modes of a
cell:\ the modes that no other mode is significantly better than, provided
that at least one mode is significantly better than another. A cell with no
bold value has no significant difference to the best. The ``Wins'' row counts the wins of
each mode, ties included. In the text we quote wins at the low, medium, and
high threshold levels as $\text{low}/\text{medium}/\text{high}$, out of $30$
configurations per level ($21$ in Tables~\ref{tab:joint-gap} and~\ref{tab:hv},
which omit the single-constraint problems). Since the tests compare the runs pairwise while the
cells show means, a few runs that stall far from the optimum can give a mode a
worse mean than a mode it beats in most runs; we return to this below.

On the best-so-far feasible objective (Table~\ref{tab:best-so-far}), even though the baseline was designed for this case, it is beaten by our alternative modes in most configurations (wins baseline 6/6/6, relaxed 13/15/16,
joint 8/10/12, ranking 19/17/19, hypervolume 12/12/10).
This supports our claim that information about the shadow-front can be obtained for free, mostly without jeopardizing (and often improving on) the performance of the baseline on the original problem.

On the four scenario metrics that aim to capture the shadow-front information, the baseline almost never wins a configuration, and an
alternative mode is significantly worse than it in only one comparison:\ the joint
mode on the whole curve under single relaxation on G24 at $N = 36$ at the high
level, by about one percent in the mean. The
relaxed, ranking, and hypervolume modes are significantly better than the
baseline in $256$, $267$, and $262$ of the $306$ comparisons with it, and
most of their remaining comparisons fall on G02 and Townsend. The joint
mode is significantly better in $240$; besides G02 and Townsend, its remaining
comparisons fall mostly on G24, G09, and the pressure vessel.

On the end points under single relaxation (Table~\ref{tab:single-gap}) the
hypervolume mode wins most often at the low and medium levels and the relaxed
mode, by one configuration, at the high level (relaxed 14/13/19, joint 12/9/11,
ranking 13/8/14, hypervolume 21/19/18). 

On the whole curve under single relaxation (Table~\ref{tab:single-hv}) the
hypervolume mode covers the per-constraint curves best (hypervolume 22/23/23,
relaxed 5/5/2, joint 3/1/0, ranking 9/9/10). On the three single-constraint problems this
metric coincides with the joint hypervolume the mode optimizes, which accounts
for eight or nine of its wins per level, but it also leads on the cantilevered
beam in every configuration, on G08 in all but one, where it ties with the
ranking mode, and on G09 in six. Since the hypervolume mode selects on a
closely related quantity, we also compare the remaining modes among themselves
(underlined:\ the best when the hypervolume mode is ignored). The ranking mode, designed for this scenario, wins most often among them at
every level (ranking 26/27/27, relaxed 6/5/2, joint 5/1/0), and it is the mode that most often wins where the hypervolume mode does not
(in 12 of 90 configurations, all nine of G24 among them).

On the end points under joint relaxation (Table~\ref{tab:joint-gap}) the
hypervolume mode leads, and the joint mode, designed for this scenario, comes
clearly second (hypervolume 17/18/16, relaxed 4/2/3, joint 9/9/10,
ranking 4/2/3). With one constraint, joint and single relaxation coincide, so
this table and Table~\ref{tab:hv} omit Branin, Mystery, and Townsend, whose
results are in Tables~\ref{tab:single-gap} and~\ref{tab:single-hv}.

On the hypervolume under joint relaxation (Table~\ref{tab:hv}) the hypervolume
mode wins every configuration outside G02 (18/18/18); on G02 the
hypervolumes are zero almost everywhere, because the runs end far from even the relaxed optimum and even with the constraints relaxed to their thresholds, find no solution as good as the strict optimum. That the hypervolume mode works well here is expected, since it selects on the metric; we therefore again 
additionally compare the remaining modes among themselves (underlined). The ranking
mode is clearly the best of these (ranking 17/17/14, relaxed 2/3/4, joint
4/3/4, baseline 0/0/0), so among the modes not driven by the metric it covers the joint
shadow-front best.

A few additional observations. On G08 the baseline is the strongest
optimizer, winning every configuration on the best-so-far objective (six tied with the relaxed mode), while
losing every scenario metric; its relative gap under single relaxation is above
$0.6$ in every configuration, against below $5 \times 10^{-3}$ for the
alternative modes, because its feasibility-first selection rarely keeps
solutions beyond the boundary of $c_2$. On Mystery above all, the means and
the tests can disagree, because a few runs of the alternative modes stall far
from the optimum and dominate the means. At the low level with $N = 36$, for
instance, the relaxed mode has a larger mean gap than the baseline
($0.67$ against $0.056$) because two of its runs stall with gaps above $1$, yet
it is better than the baseline in $95$ of the $100$ paired runs, and its median
gap is $0.002$ against $0.047$. On the pressure vessel the relative gaps of all
modes exceed $1$ at the low level under single relaxation for $N \leq 36$:\
within the budget, the modes on average do not find, even with a constraint relaxed, a solution as good as the strict optimum, so the single-constraint end points are estimated from solutions
that are still far from converged, and the modes are separated mainly by how
far they have come. Finally, on G02 the modes show significant differences at
the low level only at $N = 36$, where all modes but the hypervolume mode tie
for the win.

\begin{table}[htbp]
\centering
\caption{Final-generation mean best-so-far feasible objective by mode (Bs = Baseline, Rx = Relaxed, Jt = Joint, Rk = Ranking, Hv = Hypervolume).  \textbf{Bold}: the winning mode(s), i.e. the modes that no other mode is significantly better than, provided that some mode is significantly better than another; no mark when no mode is (paired two-sided Wilcoxon signed-rank tests, $\alpha=0.05$, Holm-corrected over the mode pairs). Lower is better. The Wins row counts the configurations each mode wins, ties included; more wins are better. On the single-constraint problems (Branin, Mystery, Townsend) the relaxed and joint modes coincide, so a single value spans their two columns; the Wins row counts it for both.}
\label{tab:best-so-far}
\small
\setlength{\tabcolsep}{2pt}
\resizebox{\ifdim\width>\linewidth\linewidth\else\width\fi}{!}{\begin{tabular}{p{0.7cm} @{\hspace{3pt}}r@{\hspace{8pt}} c c c c c @{\hspace{10pt}} c c c c c @{\hspace{10pt}} c c c c c}
\toprule
 &  & \multicolumn{5}{c}{\hrulefill\;low\;\hrulefill} & \multicolumn{5}{c}{\hrulefill\;medium\;\hrulefill} & \multicolumn{5}{c}{\hrulefill\;high\;\hrulefill} \\
 & $n$ & \rotatebox{90}{\small Bs} & \rotatebox{90}{\small Rx} & \rotatebox{90}{\small Jt} & \rotatebox{90}{\small Rk} & \rotatebox{90}{\small Hv} & \rotatebox{90}{\small Bs} & \rotatebox{90}{\small Rx} & \rotatebox{90}{\small Jt} & \rotatebox{90}{\small Rk} & \rotatebox{90}{\small Hv} & \rotatebox{90}{\small Bs} & \rotatebox{90}{\small Rx} & \rotatebox{90}{\small Jt} & \rotatebox{90}{\small Rk} & \rotatebox{90}{\small Hv} \\
\midrule
\multirow{3}{*}{\rotatebox{90}{\resizebox{\ifdim\width>1.15cm 1.15cm\else\width\fi}{!}{\shortstack[c]{\small Branin \\ \scriptsize $(\times 10^{2})$}}}} & 12 & -2.52 & \multicolumn{2}{c}{\textbf{-2.66}} & \textbf{-2.64} & -2.65 & -2.52 & \multicolumn{2}{c}{-2.67} & -2.68 & -2.69 & -2.52 & \multicolumn{2}{c}{-2.56} & \textbf{-2.69} & -2.69 \\
 & 36 & \textbf{-2.69} & \multicolumn{2}{c}{\textbf{-2.69}} & \textbf{-2.69} & -2.69 & \textbf{-2.69} & \multicolumn{2}{c}{\textbf{-2.69}} & -2.69 & -2.69 & \textbf{-2.69} & \multicolumn{2}{c}{\textbf{-2.68}} & -2.69 & -2.69 \\
 & 60 & \textbf{-2.68} & \multicolumn{2}{c}{\textbf{-2.69}} & \textbf{-2.69} & -2.69 & \textbf{-2.68} & \multicolumn{2}{c}{\textbf{-2.69}} & -2.69 & -2.69 & \textbf{-2.68} & \multicolumn{2}{c}{\textbf{-2.69}} & -2.69 & -2.69 \\
\cmidrule{2-17}
\multirow{3}{*}{\rotatebox{90}{\resizebox{\ifdim\width>1.15cm 1.15cm\else\width\fi}{!}{\small Mystery}}} & 12 & -0.542 & \multicolumn{2}{c}{-0.244} & \textbf{-0.322} & -0.298 & -0.542 & \multicolumn{2}{c}{-0.542} & -0.482 & -0.543 & -0.542 & \multicolumn{2}{c}{-0.511} & -0.585 & -0.703 \\
 & 36 & -1.17 & \multicolumn{2}{c}{-1.07} & -1.12 & -1.12 & -1.17 & \multicolumn{2}{c}{\textbf{-1.12}} & \textbf{-1.12} & -1.07 & -1.17 & \multicolumn{2}{c}{\textbf{-1.12}} & \textbf{-1.12} & -1.12 \\
 & 60 & -1.17 & \multicolumn{2}{c}{-1.17} & -1.17 & -1.17 & \textbf{-1.17} & \multicolumn{2}{c}{\textbf{-1.17}} & \textbf{-1.17} & -1.17 & \textbf{-1.17} & \multicolumn{2}{c}{\textbf{-1.17}} & \textbf{-1.17} & -1.17 \\
\cmidrule{2-17}
\multirow{3}{*}{\rotatebox{90}{\resizebox{\ifdim\width>1.15cm 1.15cm\else\width\fi}{!}{\small Townsend}}} & 12 & -1.68 & \multicolumn{2}{c}{-1.69} & -1.69 & -1.69 & -1.68 & \multicolumn{2}{c}{-1.69} & -1.69 & -1.69 & -1.68 & \multicolumn{2}{c}{\textbf{-1.71}} & \textbf{-1.72} & \textbf{-1.72} \\
 & 36 & -1.70 & \multicolumn{2}{c}{\textbf{-1.74}} & \textbf{-1.75} & \textbf{-1.75} & -1.70 & \multicolumn{2}{c}{-1.75} & -1.75 & -1.75 & -1.70 & \multicolumn{2}{c}{\textbf{-1.83}} & \textbf{-1.81} & \textbf{-1.83} \\
 & 60 & -1.74 & \multicolumn{2}{c}{\textbf{-1.77}} & \textbf{-1.77} & \textbf{-1.77} & -1.74 & \multicolumn{2}{c}{\textbf{-1.79}} & \textbf{-1.79} & \textbf{-1.79} & -1.74 & \multicolumn{2}{c}{\textbf{-1.84}} & \textbf{-1.86} & \textbf{-1.84} \\
\cmidrule{2-17}
\multirow{3}{*}{\rotatebox{90}{\resizebox{\ifdim\width>1.15cm 1.15cm\else\width\fi}{!}{\shortstack[c]{\small C.\ beam \\ \scriptsize $(\times 10^{1})$}}}} & 12 & 9.64 & 9.54 & 9.56 & 9.38 & \textbf{9.30} & 9.64 & 9.54 & 9.65 & 9.41 & \textbf{9.31} & 9.64 & 9.54 & 9.75 & 9.41 & \textbf{9.33} \\
 & 36 & 9.46 & 9.42 & 9.40 & \textbf{9.28} & \textbf{9.28} & 9.46 & 9.44 & 9.45 & \textbf{9.28} & \textbf{9.28} & 9.46 & 9.44 & 9.48 & \textbf{9.28} & 9.28 \\
 & 60 & 9.39 & 9.37 & 9.37 & \textbf{9.28} & \textbf{9.28} & 9.39 & 9.38 & 9.38 & \textbf{9.28} & 9.28 & 9.39 & 9.37 & 9.41 & \textbf{9.28} & 9.28 \\
\cmidrule{2-17}
\multirow{3}{*}{\rotatebox{90}{\resizebox{\ifdim\width>1.15cm 1.15cm\else\width\fi}{!}{\shortstack[c]{\small P.\ vessel \\ \scriptsize $(\times 10^{3})$}}}} & 12 & 6.53 & \textbf{6.48} & 6.50 & \textbf{6.45} & \textbf{6.43} & 6.53 & \textbf{6.43} & 6.51 & 6.43 & \textbf{6.37} & 6.53 & 6.46 & 6.48 & 6.36 & \textbf{6.28} \\
 & 36 & 6.59 & 6.39 & 6.43 & 6.35 & \textbf{6.24} & 6.59 & 6.35 & 6.47 & 6.29 & \textbf{6.10} & 6.59 & 6.29 & 6.41 & 6.18 & \textbf{5.99} \\
 & 60 & 6.46 & 6.32 & 6.31 & 6.27 & \textbf{6.11} & 6.46 & 6.25 & 6.31 & 6.17 & \textbf{5.97} & 6.46 & 6.14 & 6.31 & 6.06 & \textbf{5.94} \\
\cmidrule{2-17}
\multirow{3}{*}{\rotatebox{90}{\resizebox{\ifdim\width>1.15cm 1.15cm\else\width\fi}{!}{\shortstack[c]{\small G02 \\ \scriptsize $(\times 10^{-1})$}}}} & 12 & -3.90 & -3.92 & -3.92 & -3.91 & -3.89 & -3.90 & \textbf{-3.94} & \textbf{-3.95} & \textbf{-3.91} & \textbf{-3.95} & -3.90 & \textbf{-3.98} & \textbf{-4.02} & \textbf{-3.94} & \textbf{-3.99} \\
 & 36 & \textbf{-5.80} & \textbf{-5.88} & \textbf{-5.78} & \textbf{-5.83} & -5.79 & -5.80 & \textbf{-5.94} & \textbf{-5.87} & \textbf{-5.82} & \textbf{-5.88} & -5.80 & \textbf{-5.95} & \textbf{-5.93} & \textbf{-5.91} & \textbf{-5.93} \\
 & 60 & -6.60 & -6.73 & -6.68 & -6.67 & -6.64 & -6.60 & \textbf{-6.78} & \textbf{-6.78} & \textbf{-6.76} & \textbf{-6.73} & -6.60 & \textbf{-6.85} & \textbf{-6.83} & \textbf{-6.78} & \textbf{-6.74} \\
\cmidrule{2-17}
\multirow{3}{*}{\rotatebox{90}{\resizebox{\ifdim\width>1.15cm 1.15cm\else\width\fi}{!}{\shortstack[c]{\small G06 \\ \scriptsize $(\times 10^{3})$}}}} & 12 & -6.06 & \textbf{-6.42} & \textbf{-6.41} & \textbf{-6.43} & \textbf{-6.35} & -6.06 & \textbf{-6.45} & \textbf{-6.40} & \textbf{-6.51} & \textbf{-6.35} & -6.06 & \textbf{-6.51} & -6.35 & \textbf{-6.53} & -6.45 \\
 & 36 & -6.55 & \textbf{-6.75} & \textbf{-6.76} & \textbf{-6.79} & \textbf{-6.73} & -6.55 & \textbf{-6.77} & -6.69 & \textbf{-6.78} & \textbf{-6.66} & -6.55 & \textbf{-6.72} & -6.60 & \textbf{-6.76} & -6.69 \\
 & 60 & -6.71 & \textbf{-6.87} & -6.82 & \textbf{-6.86} & \textbf{-6.84} & -6.71 & -6.83 & -6.78 & \textbf{-6.88} & -6.77 & -6.71 & -6.78 & -6.74 & \textbf{-6.85} & -6.76 \\
\cmidrule{2-17}
\multirow{3}{*}{\rotatebox{90}{\resizebox{\ifdim\width>1.15cm 1.15cm\else\width\fi}{!}{\shortstack[c]{\small G08 \\ \scriptsize $(\times 10^{-2})$}}}} & 12 & \textbf{-9.58} & \textbf{-9.58} & -9.44 & -9.57 & -9.57 & \textbf{-9.58} & \textbf{-9.58} & -9.51 & -9.58 & -9.57 & \textbf{-9.58} & \textbf{-9.58} & -9.51 & -9.58 & -9.37 \\
 & 36 & \textbf{-9.58} & -9.58 & -9.58 & -9.58 & -9.58 & \textbf{-9.58} & \textbf{-9.58} & -9.58 & -9.58 & -9.57 & \textbf{-9.58} & \textbf{-9.58} & -9.58 & -9.58 & -9.57 \\
 & 60 & \textbf{-9.58} & \textbf{-9.58} & -9.58 & -9.58 & -9.58 & \textbf{-9.58} & -9.58 & -9.58 & -9.58 & -9.57 & \textbf{-9.58} & -9.58 & -9.58 & -9.58 & -9.58 \\
\cmidrule{2-17}
\multirow{3}{*}{\rotatebox{90}{\resizebox{\ifdim\width>1.15cm 1.15cm\else\width\fi}{!}{\shortstack[c]{\small G09 \\ \scriptsize $(\times 10^{2})$}}}} & 12 & 6.95 & 6.89 & 6.91 & \textbf{6.84} & \textbf{6.84} & 6.95 & 6.91 & 6.93 & \textbf{6.85} & \textbf{6.85} & 6.95 & 6.89 & 6.95 & \textbf{6.85} & 6.87 \\
 & 36 & 6.85 & 6.83 & 6.86 & \textbf{6.81} & 6.82 & 6.85 & 6.84 & 6.87 & \textbf{6.82} & 6.82 & 6.85 & 6.84 & 6.86 & \textbf{6.82} & 6.83 \\
 & 60 & 6.84 & 6.83 & 6.84 & \textbf{6.81} & 6.81 & 6.84 & 6.83 & 6.84 & \textbf{6.81} & 6.82 & 6.84 & 6.83 & 6.85 & \textbf{6.81} & 6.82 \\
\cmidrule{2-17}
\multirow{3}{*}{\rotatebox{90}{\resizebox{\ifdim\width>1.15cm 1.15cm\else\width\fi}{!}{\small G24}}} & 12 & -5.48 & -5.47 & -5.47 & \textbf{-5.50} & -5.48 & -5.48 & -5.48 & -5.49 & \textbf{-5.47} & -5.48 & -5.48 & -5.49 & -5.49 & -5.50 & -5.49 \\
 & 36 & -5.50 & -5.50 & -5.50 & \textbf{-5.51} & -5.50 & -5.50 & \textbf{-5.50} & \textbf{-5.50} & \textbf{-5.51} & -5.50 & -5.50 & \textbf{-5.50} & \textbf{-5.50} & \textbf{-5.51} & -5.50 \\
 & 60 & -5.50 & \textbf{-5.51} & -5.50 & \textbf{-5.51} & -5.51 & -5.50 & \textbf{-5.51} & -5.51 & \textbf{-5.51} & -5.50 & -5.50 & \textbf{-5.51} & \textbf{-5.51} & \textbf{-5.51} & -5.50 \\
\midrule
\multicolumn{2}{r}{\small Wins} & \textbf{6} & \textbf{13} & \textbf{8} & \textbf{19} & \textbf{12} & \textbf{6} & \textbf{15} & \textbf{10} & \textbf{17} & \textbf{12} & \textbf{6} & \textbf{16} & \textbf{12} & \textbf{19} & \textbf{10} \\
\bottomrule
\end{tabular}
}
\end{table}
 
\begin{table}[htbp]
\centering
\caption{Final-generation mean relative optimality gap of the relaxed optimum under single relaxation, averaged over the contributing constraints by mode (Bs = Baseline, Rx = Relaxed, Jt = Joint, Rk = Ranking, Hv = Hypervolume).  \textbf{Bold}: the winning mode(s), i.e. the modes that no other mode is significantly better than, provided that some mode is significantly better than another; no mark when no mode is (paired two-sided Wilcoxon signed-rank tests, $\alpha=0.05$, Holm-corrected over the mode pairs). Lower is better. Dashes mark problems without a positive reference gain. The Wins row counts the configurations each mode wins, ties included; more wins are better. On the single-constraint problems (Branin, Mystery, Townsend) the relaxed and joint modes coincide, so a single value spans their two columns; the Wins row counts it for both.}
\label{tab:single-gap}
\small
\setlength{\tabcolsep}{2pt}
\resizebox{\ifdim\width>\linewidth\linewidth\else\width\fi}{!}{\begin{tabular}{p{0.7cm} @{\hspace{3pt}}r@{\hspace{8pt}} c c c c c @{\hspace{10pt}} c c c c c @{\hspace{10pt}} c c c c c}
\toprule
 &  & \multicolumn{5}{c}{\hrulefill\;low\;\hrulefill} & \multicolumn{5}{c}{\hrulefill\;medium\;\hrulefill} & \multicolumn{5}{c}{\hrulefill\;high\;\hrulefill} \\
 & $n$ & \rotatebox{90}{\small Bs} & \rotatebox{90}{\small Rx} & \rotatebox{90}{\small Jt} & \rotatebox{90}{\small Rk} & \rotatebox{90}{\small Hv} & \rotatebox{90}{\small Bs} & \rotatebox{90}{\small Rx} & \rotatebox{90}{\small Jt} & \rotatebox{90}{\small Rk} & \rotatebox{90}{\small Hv} & \rotatebox{90}{\small Bs} & \rotatebox{90}{\small Rx} & \rotatebox{90}{\small Jt} & \rotatebox{90}{\small Rk} & \rotatebox{90}{\small Hv} \\
\midrule
\multirow{3}{*}{\rotatebox{90}{\resizebox{\ifdim\width>1.15cm 1.15cm\else\width\fi}{!}{\shortstack[c]{\small Branin \\ \scriptsize $(\times 10^{-1})$}}}} & 12 & 5.75 & \multicolumn{2}{c}{1.16} & 2.00 & \textbf{1.43} & 3.62 & \multicolumn{2}{c}{0.010} & 0.179 & \textbf{0.005} & 2.75 & \multicolumn{2}{c}{0.010} & 0.009 & \textbf{0.005} \\
 & 36 & 0.238 & \multicolumn{2}{c}{\textbf{0.004}} & 0.004 & \textbf{0.002} & 0.317 & \multicolumn{2}{c}{\textbf{0.003}} & 0.005 & \textbf{0.002} & 0.508 & \multicolumn{2}{c}{\textbf{0.003}} & 0.006 & \textbf{0.003} \\
 & 60 & 0.361 & \multicolumn{2}{c}{\textbf{0.002}} & 0.003 & \textbf{0.002} & 0.340 & \multicolumn{2}{c}{\textbf{0.001}} & 0.004 & \textbf{0.002} & 0.469 & \multicolumn{2}{c}{\textbf{0.001}} & 0.005 & 0.002 \\
\cmidrule{2-17}
\multirow{3}{*}{\rotatebox{90}{\resizebox{\ifdim\width>1.15cm 1.15cm\else\width\fi}{!}{\small Mystery}}} & 12 & 3.44 & \multicolumn{2}{c}{\textbf{5.80}} & \textbf{5.32} & \textbf{5.46} & 2.00 & \multicolumn{2}{c}{\textbf{2.62}} & 2.88 & 2.62 & 1.64 & \multicolumn{2}{c}{\textbf{2.22}} & \textbf{1.94} & 1.55 \\
 & 36 & 0.056 & \multicolumn{2}{c}{\textbf{0.671}} & \textbf{0.340} & \textbf{0.339} & 0.047 & \multicolumn{2}{c}{\textbf{0.223}} & 0.226 & 0.448 & 0.032 & \multicolumn{2}{c}{\textbf{0.180}} & 0.182 & 0.183 \\
 & 60 & 0.040 & \multicolumn{2}{c}{\textbf{0.003}} & 0.006 & 0.004 & 0.028 & \multicolumn{2}{c}{\textbf{0.001}} & 0.003 & 0.003 & 0.024 & \multicolumn{2}{c}{\textbf{0.000}} & 0.001 & 0.002 \\
\cmidrule{2-17}
\multirow{3}{*}{\rotatebox{90}{\resizebox{\ifdim\width>1.15cm 1.15cm\else\width\fi}{!}{\small Townsend}}} & 12 & 4.34 & \multicolumn{2}{c}{4.27} & 4.28 & 4.27 & 2.65 & \multicolumn{2}{c}{2.56} & 2.57 & 2.57 & 1.78 & \multicolumn{2}{c}{\textbf{1.64}} & \textbf{1.62} & \textbf{1.64} \\
 & 36 & 3.96 & \multicolumn{2}{c}{\textbf{3.59}} & \textbf{3.45} & \textbf{3.54} & 2.37 & \multicolumn{2}{c}{2.15} & 2.15 & 2.12 & 1.53 & \multicolumn{2}{c}{\textbf{1.04}} & \textbf{1.14} & \textbf{1.06} \\
 & 60 & 3.53 & \multicolumn{2}{c}{\textbf{3.29}} & \textbf{3.20} & \textbf{3.20} & 2.08 & \multicolumn{2}{c}{1.82} & 1.80 & \textbf{1.80} & 1.39 & \multicolumn{2}{c}{\textbf{0.962}} & \textbf{0.902} & \textbf{0.981} \\
\cmidrule{2-17}
\multirow{3}{*}{\rotatebox{90}{\resizebox{\ifdim\width>1.15cm 1.15cm\else\width\fi}{!}{\shortstack[c]{\small C.\ beam \\ \scriptsize $(\times 10^{-1})$}}}} & 12 & 1.33 & 0.342 & 0.405 & 0.136 & \textbf{0.028} & 1.73 & 0.248 & 0.714 & 0.100 & \textbf{0.030} & 1.73 & 0.248 & 1.25 & 0.100 & \textbf{0.024} \\
 & 36 & 0.951 & 0.135 & 0.129 & 0.015 & \textbf{0.003} & 1.16 & 0.126 & 0.168 & 0.019 & \textbf{0.005} & 1.16 & 0.126 & 0.282 & 0.019 & \textbf{0.006} \\
 & 60 & 0.788 & 0.091 & 0.076 & 0.017 & \textbf{0.001} & 0.990 & 0.076 & 0.101 & 0.018 & \textbf{0.003} & 0.990 & 0.075 & 0.197 & 0.019 & \textbf{0.003} \\
\cmidrule{2-17}
\multirow{3}{*}{\rotatebox{90}{\resizebox{\ifdim\width>1.15cm 1.15cm\else\width\fi}{!}{\small P.\ vessel}}} & 12 & 1.61 & \textbf{1.44} & 1.53 & \textbf{1.39} & \textbf{1.46} & 0.768 & \textbf{0.614} & 0.739 & \textbf{0.628} & \textbf{0.672} & 0.447 & \textbf{0.334} & 0.410 & \textbf{0.296} & 0.356 \\
 & 36 & 1.71 & \textbf{1.23} & 1.34 & \textbf{1.14} & \textbf{1.02} & 0.798 & 0.531 & 0.659 & 0.477 & \textbf{0.360} & 0.432 & \textbf{0.247} & 0.332 & \textbf{0.222} & \textbf{0.208} \\
 & 60 & 1.39 & 1.06 & 1.04 & 0.965 & \textbf{0.662} & 0.650 & 0.423 & 0.479 & 0.352 & \textbf{0.201} & 0.364 & \textbf{0.189} & 0.278 & \textbf{0.167} & 0.182 \\
\cmidrule{2-17}
\multirow{3}{*}{\rotatebox{90}{\resizebox{\ifdim\width>1.15cm 1.15cm\else\width\fi}{!}{\shortstack[c]{\small G02 \\ \scriptsize $(\times 10^{2})$}}}} & 12 & 1.72 & 1.71 & 1.71 & 1.71 & 1.72 & 0.835 & 0.828 & 0.827 & 0.834 & 0.827 & 0.541 & \textbf{0.531} & \textbf{0.526} & \textbf{0.536} & \textbf{0.531} \\
 & 36 & \textbf{0.930} & \textbf{0.897} & \textbf{0.938} & \textbf{0.919} & 0.934 & 0.453 & 0.427 & 0.440 & 0.451 & 0.438 & 0.294 & 0.276 & 0.279 & 0.280 & 0.278 \\
 & 60 & 0.595 & 0.544 & 0.565 & 0.567 & 0.579 & 0.290 & \textbf{0.256} & \textbf{0.257} & \textbf{0.259} & \textbf{0.267} & 0.188 & \textbf{0.158} & \textbf{0.160} & \textbf{0.166} & \textbf{0.172} \\
\cmidrule{2-17}
\multirow{3}{*}{\rotatebox{90}{\resizebox{\ifdim\width>1.15cm 1.15cm\else\width\fi}{!}{\small G06}}} & 12 & 4.27 & \textbf{2.91} & \textbf{2.84} & \textbf{2.75} & \textbf{3.24} & 2.39 & \textbf{1.51} & \textbf{1.36} & \textbf{1.42} & \textbf{1.47} & 1.68 & \textbf{0.888} & \textbf{0.923} & \textbf{0.915} & \textbf{0.808} \\
 & 36 & 1.85 & 1.18 & \textbf{1.01} & \textbf{1.03} & \textbf{1.04} & 1.08 & \textbf{0.581} & \textbf{0.544} & \textbf{0.512} & \textbf{0.518} & 0.762 & 0.375 & 0.451 & 0.352 & \textbf{0.314} \\
 & 60 & 1.30 & 0.754 & \textbf{0.697} & \textbf{0.666} & \textbf{0.604} & 0.723 & 0.367 & 0.426 & 0.378 & \textbf{0.279} & 0.511 & 0.260 & 0.376 & 0.240 & \textbf{0.157} \\
\cmidrule{2-17}
\multirow{3}{*}{\rotatebox{90}{\resizebox{\ifdim\width>1.15cm 1.15cm\else\width\fi}{!}{\shortstack[c]{\small G08 \\ \scriptsize $(\times 10^{-1})$}}}} & 12 & 6.84 & 0.008 & 0.001 & 0.047 & \textbf{0.000} & 6.84 & 0.028 & 0.001 & 0.002 & \textbf{0.000} & 6.84 & 0.005 & 0.034 & 0.047 & \textbf{0.000} \\
 & 36 & 6.64 & 0.000 & 0.000 & 0.001 & \textbf{0.000} & 6.64 & 0.000 & 0.000 & 0.001 & \textbf{0.000} & 6.64 & 0.000 & 0.006 & 0.001 & \textbf{0.000} \\
 & 60 & 6.08 & 0.000 & \textbf{0.000} & 0.001 & \textbf{0.000} & 6.08 & 0.000 & \textbf{0.000} & 0.001 & \textbf{0.000} & 6.08 & 0.000 & 0.003 & 0.001 & \textbf{0.000} \\
\cmidrule{2-17}
\multirow{3}{*}{\rotatebox{90}{\resizebox{\ifdim\width>1.15cm 1.15cm\else\width\fi}{!}{\small G09}}} & 12 & 3.73 & 2.27 & 2.77 & 1.26 & \textbf{1.02} & 3.61 & 2.78 & 3.08 & \textbf{1.39} & \textbf{1.35} & 3.62 & 2.43 & 3.65 & \textbf{1.49} & \textbf{1.64} \\
 & 36 & 1.46 & 0.924 & 1.43 & \textbf{0.589} & \textbf{0.635} & 1.39 & 0.870 & 1.73 & \textbf{0.645} & 0.740 & 1.40 & \textbf{0.949} & 1.70 & \textbf{0.719} & 0.802 \\
 & 60 & 1.14 & 0.730 & 1.09 & \textbf{0.455} & 0.522 & 1.09 & \textbf{0.628} & 1.14 & \textbf{0.527} & \textbf{0.573} & 1.09 & \textbf{0.591} & 1.27 & \textbf{0.563} & 0.676 \\
\cmidrule{2-17}
\multirow{3}{*}{\rotatebox{90}{\resizebox{\ifdim\width>1.15cm 1.15cm\else\width\fi}{!}{\small G24}}} & 12 & 0.878 & \textbf{0.907} & 1.15 & \textbf{0.263} & 0.674 & 0.513 & \textbf{0.351} & 0.406 & \textbf{0.449} & 0.396 & 0.329 & \textbf{0.094} & 0.264 & \textbf{0.088} & 0.166 \\
 & 36 & 0.383 & \textbf{0.092} & 0.341 & 0.121 & 0.180 & 0.220 & \textbf{0.051} & 0.232 & 0.074 & 0.124 & 0.174 & \textbf{0.021} & 0.144 & 0.047 & 0.089 \\
 & 60 & 0.279 & \textbf{0.048} & 0.276 & 0.082 & 0.133 & 0.182 & \textbf{0.018} & 0.192 & 0.050 & 0.092 & 0.128 & \textbf{0.008} & 0.120 & 0.033 & 0.065 \\
\midrule
\multicolumn{2}{r}{\small Wins} & \textbf{1} & \textbf{14} & \textbf{12} & \textbf{13} & \textbf{21} & \textbf{0} & \textbf{13} & \textbf{9} & \textbf{8} & \textbf{19} & \textbf{0} & \textbf{19} & \textbf{11} & \textbf{14} & \textbf{18} \\
\bottomrule
\end{tabular}
}
\end{table}
 
\begin{table}[htbp]
\centering
\caption{Final-generation mean hypervolume of the shadow-front under single relaxation, averaged over the contributing constraints by mode (Bs = Baseline, Rx = Relaxed, Jt = Joint, Rk = Ranking, Hv = Hypervolume).  \textbf{Bold}: the winning mode(s), i.e. the modes that no other mode is significantly better than, provided that some mode is significantly better than another; no mark when no mode is (paired two-sided Wilcoxon signed-rank tests, $\alpha=0.05$, Holm-corrected over the mode pairs). Higher is better. Dashes mark problems without a positive reference gain. The Wins row counts the configurations each mode wins, ties included; more wins are better. The hypervolume (Hv) mode optimizes this metric by construction;\ \underline{underline} marks the winning mode(s) as if the hypervolume mode had not been run, by the same rule on the pairs of the other modes (Holm-corrected over those pairs), and the Wins row gives these counts in parentheses next to each mode's win count (ties included). On the single-constraint problems (Branin, Mystery, Townsend) the relaxed and joint modes coincide, so a single value spans their two columns; the Wins row counts it for both.}
\label{tab:single-hv}
\small
\setlength{\tabcolsep}{2pt}
\resizebox{\ifdim\width>\linewidth\linewidth\else\width\fi}{!}{\begin{tabular}{p{0.7cm} @{\hspace{3pt}}r@{\hspace{8pt}} c c c c c @{\hspace{10pt}} c c c c c @{\hspace{10pt}} c c c c c}
\toprule
 &  & \multicolumn{5}{c}{\hrulefill\;low\;\hrulefill} & \multicolumn{5}{c}{\hrulefill\;medium\;\hrulefill} & \multicolumn{5}{c}{\hrulefill\;high\;\hrulefill} \\
 & $n$ & \rotatebox{90}{\small Bs} & \rotatebox{90}{\small Rx} & \rotatebox{90}{\small Jt} & \rotatebox{90}{\small Rk} & \rotatebox{90}{\small Hv} & \rotatebox{90}{\small Bs} & \rotatebox{90}{\small Rx} & \rotatebox{90}{\small Jt} & \rotatebox{90}{\small Rk} & \rotatebox{90}{\small Hv} & \rotatebox{90}{\small Bs} & \rotatebox{90}{\small Rx} & \rotatebox{90}{\small Jt} & \rotatebox{90}{\small Rk} & \rotatebox{90}{\small Hv} \\
\midrule
\multirow{3}{*}{\rotatebox{90}{\resizebox{\ifdim\width>1.15cm 1.15cm\else\width\fi}{!}{\shortstack[c]{\small Branin \\ \scriptsize $(\times 10^{-1})$}}}} & 12 & 4.51 & \multicolumn{2}{c}{6.16} & \underline{5.98} & \textbf{6.13} & 4.29 & \multicolumn{2}{c}{6.02} & \underline{5.99} & \textbf{6.07} & 4.08 & \multicolumn{2}{c}{5.56} & \underline{5.66} & \textbf{5.68} \\
 & 36 & 6.29 & \multicolumn{2}{c}{6.45} & \underline{6.46} & \textbf{6.46} & 5.92 & \multicolumn{2}{c}{6.06} & \underline{6.08} & \textbf{6.08} & 5.51 & \multicolumn{2}{c}{5.67} & \underline{5.69} & \textbf{5.69} \\
 & 60 & 6.27 & \multicolumn{2}{c}{6.46} & \underline{6.46} & \textbf{6.47} & 5.90 & \multicolumn{2}{c}{6.07} & \underline{6.08} & \textbf{6.08} & 5.50 & \multicolumn{2}{c}{5.68} & \underline{5.69} & \textbf{5.69} \\
\cmidrule{2-17}
\multirow{3}{*}{\rotatebox{90}{\resizebox{\ifdim\width>1.15cm 1.15cm\else\width\fi}{!}{\shortstack[c]{\small Mystery \\ \scriptsize $(\times 10^{-1})$}}}} & 12 & 4.24 & \multicolumn{2}{c}{4.34} & \underline{4.63} & \textbf{4.63} & 4.64 & \multicolumn{2}{c}{4.86} & \underline{5.07} & \textbf{5.11} & 5.05 & \multicolumn{2}{c}{5.22} & \underline{5.55} & \textbf{5.65} \\
 & 36 & 5.11 & \multicolumn{2}{c}{5.17} & \underline{5.40} & \textbf{5.43} & 5.46 & \multicolumn{2}{c}{5.52} & \underline{5.69} & \textbf{5.65} & 5.90 & \multicolumn{2}{c}{5.92} & \underline{6.09} & \textbf{6.10} \\
 & 60 & 5.22 & \multicolumn{2}{c}{5.33} & \underline{5.48} & \textbf{5.50} & 5.56 & \multicolumn{2}{c}{5.64} & \underline{5.77} & \textbf{5.78} & 5.98 & \multicolumn{2}{c}{6.04} & \underline{6.17} & \textbf{6.18} \\
\cmidrule{2-17}
\multirow{3}{*}{\rotatebox{90}{\resizebox{\ifdim\width>1.15cm 1.15cm\else\width\fi}{!}{\shortstack[c]{\small Townsend \\ \scriptsize $(\times 10^{-1})$}}}} & 12 & 0.301 & \multicolumn{2}{c}{0.439} & 0.433 & 0.454 & 0.324 & \multicolumn{2}{c}{0.561} & \underline{0.566} & \textbf{0.571} & 0.380 & \multicolumn{2}{c}{0.883} & \underline{0.944} & \textbf{0.897} \\
 & 36 & 0.574 & \multicolumn{2}{c}{1.18} & \underline{1.35} & \textbf{1.25} & 0.632 & \multicolumn{2}{c}{1.29} & \underline{1.31} & \textbf{1.36} & 0.719 & \multicolumn{2}{c}{2.48} & \underline{2.23} & \textbf{2.44} \\
 & 60 & 1.12 & \multicolumn{2}{c}{1.49} & \underline{1.61} & \textbf{1.62} & 1.15 & \multicolumn{2}{c}{1.87} & \underline{1.94} & \textbf{1.94} & 1.22 & \multicolumn{2}{c}{2.69} & \underline{2.87} & \textbf{2.66} \\
\cmidrule{2-17}
\multirow{3}{*}{\rotatebox{90}{\resizebox{\ifdim\width>1.15cm 1.15cm\else\width\fi}{!}{\shortstack[c]{\small C.\ beam \\ \scriptsize $(\times 10^{-1})$}}}} & 12 & 4.83 & 5.31 & 5.24 & \underline{5.64} & \textbf{5.84} & 6.30 & 7.22 & 6.82 & \underline{7.42} & \textbf{7.55} & 7.02 & 8.13 & 7.25 & \underline{8.31} & \textbf{8.42} \\
 & 36 & 5.32 & 5.62 & 5.64 & \underline{5.86} & \textbf{5.87} & 6.79 & 7.40 & 7.34 & \underline{7.58} & \textbf{7.59} & 7.54 & 8.29 & 8.12 & \underline{8.44} & \textbf{8.46} \\
 & 60 & 5.51 & 5.71 & 5.73 & \underline{5.86} & \textbf{5.87} & 6.96 & 7.48 & 7.44 & \underline{7.58} & \textbf{7.59} & 7.70 & 8.36 & 8.24 & \underline{8.45} & \textbf{8.46} \\
\cmidrule{2-17}
\multirow{3}{*}{\rotatebox{90}{\resizebox{\ifdim\width>1.15cm 1.15cm\else\width\fi}{!}{\shortstack[c]{\small P.\ vessel \\ \scriptsize $(\times 10^{-1})$}}}} & 12 & 0.531 & \underline{\textbf{0.648}} & 0.537 & \underline{\textbf{0.762}} & \textbf{0.736} & 1.26 & \underline{\textbf{1.55}} & 1.31 & \underline{\textbf{1.56}} & \textbf{1.60} & 2.32 & 2.51 & 2.33 & \underline{\textbf{2.78}} & \textbf{2.78} \\
 & 36 & 0.440 & \underline{1.03} & \underline{0.796} & \underline{1.02} & \textbf{1.25} & 1.16 & 1.91 & 1.63 & \underline{2.21} & \textbf{2.95} & 2.29 & 3.13 & 2.74 & \underline{3.47} & \textbf{4.07} \\
 & 60 & 0.856 & \underline{1.10} & \underline{1.06} & \underline{1.43} & \textbf{2.06} & 1.67 & 2.33 & 2.10 & \underline{2.77} & \textbf{3.93} & 2.69 & 3.66 & 3.03 & \underline{3.97} & \textbf{4.34} \\
\cmidrule{2-17}
\multirow{3}{*}{\rotatebox{90}{\resizebox{\ifdim\width>1.15cm 1.15cm\else\width\fi}{!}{\shortstack[c]{\small G02 \\ \scriptsize $(\times 10^{-4})$}}}} & 12 & 0 & 0 & 0 & 0 & 0 & 0 & 0 & 0 & 0 & 0 & 0 & 0 & 0 & 0 & 0 \\
 & 36 & 0 & 0 & 0 & 0 & 0 & 0 & 0 & 0 & 0 & 0 & 0 & 0 & 0 & 0 & 0 \\
 & 60 & 0 & 0 & 0 & 0 & 0 & 0 & 2.80 & 0 & 0 & 0 & 0 & 0 & 0 & 0 & 0 \\
\cmidrule{2-17}
\multirow{3}{*}{\rotatebox{90}{\resizebox{\ifdim\width>1.15cm 1.15cm\else\width\fi}{!}{\shortstack[c]{\small G06 \\ \scriptsize $(\times 10^{-1})$}}}} & 12 & 0.080 & \underline{\textbf{0.235}} & \underline{\textbf{0.292}} & \underline{\textbf{0.351}} & \textbf{0.304} & 0.199 & \underline{\textbf{0.786}} & \underline{\textbf{0.983}} & \underline{\textbf{0.730}} & \textbf{0.912} & 0.364 & \underline{\textbf{1.51}} & 1.22 & \underline{\textbf{1.58}} & \textbf{1.77} \\
 & 36 & 0.255 & \underline{\textbf{1.17}} & \underline{\textbf{1.21}} & \underline{\textbf{1.13}} & \textbf{1.49} & 0.674 & \underline{\textbf{2.29}} & 2.03 & \underline{\textbf{2.49}} & \textbf{2.53} & 1.15 & 3.01 & 2.12 & \underline{\textbf{3.31}} & \textbf{3.13} \\
 & 60 & 0.364 & \underline{\textbf{1.69}} & \underline{\textbf{1.69}} & \underline{\textbf{2.05}} & \textbf{1.97} & 1.12 & \underline{\textbf{3.12}} & 2.36 & \underline{\textbf{3.26}} & \textbf{3.32} & 1.77 & 3.50 & 2.44 & \underline{\textbf{3.87}} & \textbf{3.64} \\
\cmidrule{2-17}
\multirow{3}{*}{\rotatebox{90}{\resizebox{\ifdim\width>1.15cm 1.15cm\else\width\fi}{!}{\shortstack[c]{\small G08 \\ \scriptsize $(\times 10^{-1})$}}}} & 12 & 2.89 & 5.46 & 5.49 & \underline{5.51} & \textbf{5.53} & 3.02 & 7.03 & 7.06 & \underline{7.08} & \textbf{7.09} & 3.09 & 7.83 & 7.77 & \underline{\textbf{7.83}} & \textbf{7.86} \\
 & 36 & 3.04 & 5.52 & 5.52 & \underline{5.54} & \textbf{5.54} & 3.19 & 7.08 & 7.08 & \underline{7.09} & \textbf{7.10} & 3.27 & 7.86 & 7.83 & \underline{7.87} & \textbf{7.87} \\
 & 60 & 3.27 & 5.53 & 5.53 & \underline{5.54} & \textbf{5.54} & 3.53 & 7.09 & 7.09 & \underline{7.09} & \textbf{7.10} & 3.67 & 7.86 & 7.85 & \underline{7.87} & \textbf{7.87} \\
\cmidrule{2-17}
\multirow{3}{*}{\rotatebox{90}{\resizebox{\ifdim\width>1.15cm 1.15cm\else\width\fi}{!}{\shortstack[c]{\small G09 \\ \scriptsize $(\times 10^{-1})$}}}} & 12 & 2.20 & 2.39 & 2.32 & \underline{2.59} & \textbf{2.70} & 2.31 & 2.55 & 2.51 & \underline{2.78} & \textbf{2.85} & 2.30 & 2.72 & 2.57 & \underline{2.93} & \textbf{2.98} \\
 & 36 & 2.46 & \textbf{2.90} & 2.63 & \underline{\textbf{3.02}} & \textbf{2.84} & 2.56 & 3.02 & 2.73 & \underline{\textbf{3.21}} & \textbf{2.97} & 2.60 & 3.13 & 2.82 & \underline{\textbf{3.39}} & \textbf{3.09} \\
 & 60 & 2.54 & 3.07 & 2.68 & \underline{\textbf{3.33}} & 3.03 & 2.65 & \underline{\textbf{3.51}} & 2.86 & \underline{\textbf{3.44}} & 3.10 & 2.71 & \underline{\textbf{3.70}} & 2.92 & \underline{\textbf{3.55}} & 3.19 \\
\cmidrule{2-17}
\multirow{3}{*}{\rotatebox{90}{\resizebox{\ifdim\width>1.15cm 1.15cm\else\width\fi}{!}{\shortstack[c]{\small G24 \\ \scriptsize $(\times 10^{-1})$}}}} & 12 & 2.27 & 3.11 & 2.42 & \underline{\textbf{3.58}} & 3.16 & 2.70 & 3.59 & 2.82 & \underline{\textbf{3.97}} & 3.51 & 3.30 & 4.02 & 3.26 & \underline{\textbf{4.39}} & 3.93 \\
 & 36 & 2.85 & 3.89 & 3.10 & \underline{\textbf{4.50}} & 4.19 & 3.52 & 4.25 & 3.59 & \underline{\textbf{4.70}} & 4.41 & 4.05 & 4.40 & 4.00 & \underline{\textbf{4.84}} & 4.58 \\
 & 60 & 3.26 & 4.20 & 3.35 & \underline{\textbf{4.75}} & 4.50 & 3.89 & 4.36 & 3.87 & \underline{\textbf{4.89}} & 4.64 & 4.29 & 4.57 & 4.26 & \underline{\textbf{4.96}} & 4.74 \\
\midrule
\multicolumn{2}{r}{\small Wins} & \textbf{0} (0) & \textbf{5} (6) & \textbf{3} (5) & \textbf{9} (26) & \textbf{22} & \textbf{0} (0) & \textbf{5} (5) & \textbf{1} (1) & \textbf{9} (27) & \textbf{23} & \textbf{0} (0) & \textbf{2} (2) & \textbf{0} (0) & \textbf{10} (27) & \textbf{23} \\
\bottomrule
\end{tabular}
}
\end{table}
 
\begin{table}[htbp]
\centering
\caption{Final-generation mean relative optimality gap of the threshold optimum under joint relaxation by mode (Bs = Baseline, Rx = Relaxed, Jt = Joint, Rk = Ranking, Hv = Hypervolume).  \textbf{Bold}: the winning mode(s), i.e. the modes that no other mode is significantly better than, provided that some mode is significantly better than another; no mark when no mode is (paired two-sided Wilcoxon signed-rank tests, $\alpha=0.05$, Holm-corrected over the mode pairs). Lower is better. Dashes mark problems without a positive reference gain. The Wins row counts the configurations each mode wins, ties included; more wins are better. The single-constraint problems (Branin, Mystery, Townsend) are omitted: with one constraint, joint and single relaxation coincide, and their results are in Table~\ref{tab:single-gap}.}
\label{tab:joint-gap}
\small
\setlength{\tabcolsep}{2pt}
\resizebox{\ifdim\width>\linewidth\linewidth\else\width\fi}{!}{\begin{tabular}{p{0.7cm} @{\hspace{3pt}}r@{\hspace{8pt}} c c c c c @{\hspace{10pt}} c c c c c @{\hspace{10pt}} c c c c c}
\toprule
 &  & \multicolumn{5}{c}{\hrulefill\;low\;\hrulefill} & \multicolumn{5}{c}{\hrulefill\;medium\;\hrulefill} & \multicolumn{5}{c}{\hrulefill\;high\;\hrulefill} \\
 & $n$ & \rotatebox{90}{\small Bs} & \rotatebox{90}{\small Rx} & \rotatebox{90}{\small Jt} & \rotatebox{90}{\small Rk} & \rotatebox{90}{\small Hv} & \rotatebox{90}{\small Bs} & \rotatebox{90}{\small Rx} & \rotatebox{90}{\small Jt} & \rotatebox{90}{\small Rk} & \rotatebox{90}{\small Hv} & \rotatebox{90}{\small Bs} & \rotatebox{90}{\small Rx} & \rotatebox{90}{\small Jt} & \rotatebox{90}{\small Rk} & \rotatebox{90}{\small Hv} \\
\midrule
\multirow{3}{*}{\rotatebox{90}{\resizebox{\ifdim\width>1.15cm 1.15cm\else\width\fi}{!}{\shortstack[c]{\small C.\ beam \\ \scriptsize $(\times 10^{-1})$}}}} & 12 & 1.33 & 0.342 & 0.403 & 0.136 & \textbf{0.028} & 2.22 & 0.262 & 0.204 & 0.180 & \textbf{0.021} & 3.06 & 0.616 & 0.180 & 0.630 & \textbf{0.019} \\
 & 36 & 0.951 & 0.135 & 0.129 & 0.015 & \textbf{0.003} & 1.54 & 0.075 & 0.075 & 0.141 & \textbf{0.003} & 1.99 & 0.257 & 0.030 & 0.506 & \textbf{0.003} \\
 & 60 & 0.788 & 0.091 & 0.076 & 0.017 & \textbf{0.001} & 1.28 & 0.054 & 0.042 & 0.135 & \textbf{0.002} & 1.56 & 0.183 & 0.017 & 0.503 & \textbf{0.001} \\
\cmidrule{2-17}
\multirow{3}{*}{\rotatebox{90}{\resizebox{\ifdim\width>1.15cm 1.15cm\else\width\fi}{!}{\shortstack[c]{\small P.\ vessel \\ \scriptsize $(\times 10^{-1})$}}}} & 12 & 7.72 & 5.18 & \textbf{4.42} & 5.27 & \textbf{4.54} & 5.40 & 3.37 & \textbf{2.73} & 3.49 & \textbf{2.66} & 3.92 & 1.49 & \textbf{0.760} & 1.32 & \textbf{0.647} \\
 & 36 & 7.52 & 4.29 & \textbf{3.92} & 4.67 & \textbf{3.41} & 5.28 & 2.82 & 2.52 & 3.06 & \textbf{1.49} & 3.43 & 0.991 & 0.541 & 0.988 & \textbf{0.038} \\
 & 60 & 6.37 & 3.75 & 3.13 & 4.13 & \textbf{2.33} & 4.50 & 2.37 & 1.92 & 2.58 & \textbf{0.746} & 2.73 & 0.595 & 0.316 & 0.638 & \textbf{0.026} \\
\cmidrule{2-17}
\multirow{3}{*}{\rotatebox{90}{\resizebox{\ifdim\width>1.15cm 1.15cm\else\width\fi}{!}{\shortstack[c]{\small G02 \\ \scriptsize $(\times 10^{2})$}}}} & 12 & 1.72 & 1.71 & 1.71 & 1.71 & 1.72 & 0.835 & 0.828 & 0.827 & 0.834 & 0.827 & 0.541 & \textbf{0.531} & \textbf{0.526} & \textbf{0.536} & \textbf{0.531} \\
 & 36 & \textbf{0.930} & \textbf{0.897} & \textbf{0.938} & \textbf{0.919} & 0.934 & 0.453 & 0.427 & 0.440 & 0.451 & 0.438 & 0.294 & 0.276 & 0.279 & 0.280 & 0.278 \\
 & 60 & 0.595 & 0.544 & 0.565 & 0.567 & 0.579 & 0.290 & \textbf{0.256} & \textbf{0.257} & \textbf{0.259} & \textbf{0.267} & 0.188 & \textbf{0.158} & \textbf{0.160} & \textbf{0.166} & \textbf{0.172} \\
\cmidrule{2-17}
\multirow{3}{*}{\rotatebox{90}{\resizebox{\ifdim\width>1.15cm 1.15cm\else\width\fi}{!}{\small G06}}} & 12 & 2.04 & \textbf{1.71} & \textbf{1.65} & \textbf{1.60} & 1.99 & 1.19 & \textbf{0.865} & \textbf{0.790} & \textbf{0.868} & \textbf{0.973} & 0.844 & \textbf{0.463} & \textbf{0.426} & \textbf{0.532} & \textbf{0.408} \\
 & 36 & 0.896 & \textbf{0.782} & \textbf{0.758} & \textbf{0.700} & \textbf{0.704} & 0.529 & 0.376 & \textbf{0.307} & 0.359 & \textbf{0.288} & 0.300 & 0.142 & \textbf{0.056} & 0.155 & 0.070 \\
 & 60 & 0.633 & \textbf{0.471} & \textbf{0.406} & \textbf{0.478} & \textbf{0.460} & 0.383 & 0.249 & \textbf{0.188} & 0.236 & \textbf{0.199} & 0.190 & 0.084 & \textbf{0.013} & 0.086 & 0.033 \\
\cmidrule{2-17}
\multirow{3}{*}{\rotatebox{90}{\resizebox{\ifdim\width>1.15cm 1.15cm\else\width\fi}{!}{\shortstack[c]{\small G08 \\ \scriptsize $(\times 10^{-1})$}}}} & 12 & 6.84 & 0.008 & 0.001 & 0.047 & \textbf{0.000} & 6.84 & 0.028 & 0.001 & 0.002 & \textbf{0.000} & 9.36 & 2.22 & 0.069 & 2.16 & \textbf{0.046} \\
 & 36 & 6.64 & 0.000 & 0.000 & 0.001 & \textbf{0.000} & 6.64 & 0.000 & 0.000 & 0.001 & \textbf{0.000} & 9.13 & 0.543 & 0.015 & 1.43 & \textbf{0.005} \\
 & 60 & 6.08 & 0.000 & \textbf{0.000} & 0.001 & \textbf{0.000} & 6.08 & 0.000 & \textbf{0.000} & 0.001 & \textbf{0.000} & 9.18 & 0.293 & 0.009 & 1.58 & \textbf{0.004} \\
\cmidrule{2-17}
\multirow{3}{*}{\rotatebox{90}{\resizebox{\ifdim\width>1.15cm 1.15cm\else\width\fi}{!}{\shortstack[c]{\small G09 \\ \scriptsize $(\times 10^{-1})$}}}} & 12 & 2.11 & 1.22 & 0.856 & 0.807 & \textbf{0.488} & 2.85 & 1.24 & 0.816 & 0.730 & \textbf{0.441} & 3.38 & 1.40 & 0.799 & 0.732 & \textbf{0.355} \\
 & 36 & 1.73 & 0.826 & 0.445 & 0.556 & \textbf{0.201} & 2.27 & 0.793 & 0.393 & 0.540 & \textbf{0.147} & 2.62 & 0.707 & 0.340 & 0.525 & \textbf{0.141} \\
 & 60 & 1.54 & 0.667 & 0.351 & 0.504 & \textbf{0.127} & 1.88 & 0.597 & 0.289 & 0.482 & \textbf{0.112} & 2.16 & 0.601 & \textbf{0.160} & 0.489 & \textbf{0.108} \\
\cmidrule{2-17}
\multirow{3}{*}{\rotatebox{90}{\resizebox{\ifdim\width>1.15cm 1.15cm\else\width\fi}{!}{\shortstack[c]{\small G24 \\ \scriptsize $(\times 10^{-1})$}}}} & 12 & 2.33 & 3.21 & 2.46 & 1.12 & \textbf{1.53} & 1.76 & 1.68 & \textbf{0.435} & 1.66 & \textbf{0.778} & 2.19 & 0.973 & \textbf{0.202} & 0.778 & \textbf{0.221} \\
 & 36 & 1.14 & 0.948 & \textbf{0.264} & 0.436 & \textbf{0.175} & 1.03 & 0.724 & \textbf{0.145} & 0.420 & \textbf{0.140} & 1.28 & 0.577 & \textbf{0.055} & 0.401 & 0.118 \\
 & 60 & 0.939 & 0.892 & \textbf{0.145} & 0.326 & \textbf{0.125} & 0.770 & 0.572 & \textbf{0.055} & 0.275 & 0.116 & 1.13 & 0.490 & \textbf{0.025} & 0.321 & 0.083 \\
\midrule
\multicolumn{2}{r}{\small Wins} & \textbf{1} & \textbf{4} & \textbf{9} & \textbf{4} & \textbf{17} & \textbf{0} & \textbf{2} & \textbf{9} & \textbf{2} & \textbf{18} & \textbf{0} & \textbf{3} & \textbf{10} & \textbf{3} & \textbf{16} \\
\bottomrule
\end{tabular}
}
\end{table}
 
\begin{table}[htbp]
\centering
\caption{Final-generation mean hypervolume of the shadow-front under joint relaxation by mode (Bs = Baseline, Rx = Relaxed, Jt = Joint, Rk = Ranking, Hv = Hypervolume).  \textbf{Bold}: the winning mode(s), i.e. the modes that no other mode is significantly better than, provided that some mode is significantly better than another; no mark when no mode is (paired two-sided Wilcoxon signed-rank tests, $\alpha=0.05$, Holm-corrected over the mode pairs). Higher is better. The Wins row counts the configurations each mode wins, ties included; more wins are better. The hypervolume (Hv) mode optimizes this metric by construction;\ \underline{underline} marks the winning mode(s) as if the hypervolume mode had not been run, by the same rule on the pairs of the other modes (Holm-corrected over those pairs), and the Wins row gives these counts in parentheses next to each mode's win count (ties included). The single-constraint problems (Branin, Mystery, Townsend) are omitted: with one constraint, joint and single relaxation coincide, and their results are in Table~\ref{tab:single-hv}.}
\label{tab:hv}
\small
\setlength{\tabcolsep}{2pt}
\resizebox{\ifdim\width>\linewidth\linewidth\else\width\fi}{!}{\begin{tabular}{p{0.7cm} @{\hspace{3pt}}r@{\hspace{8pt}} c c c c c @{\hspace{10pt}} c c c c c @{\hspace{10pt}} c c c c c}
\toprule
 &  & \multicolumn{5}{c}{\hrulefill\;low\;\hrulefill} & \multicolumn{5}{c}{\hrulefill\;medium\;\hrulefill} & \multicolumn{5}{c}{\hrulefill\;high\;\hrulefill} \\
 & $n$ & \rotatebox{90}{\small Bs} & \rotatebox{90}{\small Rx} & \rotatebox{90}{\small Jt} & \rotatebox{90}{\small Rk} & \rotatebox{90}{\small Hv} & \rotatebox{90}{\small Bs} & \rotatebox{90}{\small Rx} & \rotatebox{90}{\small Jt} & \rotatebox{90}{\small Rk} & \rotatebox{90}{\small Hv} & \rotatebox{90}{\small Bs} & \rotatebox{90}{\small Rx} & \rotatebox{90}{\small Jt} & \rotatebox{90}{\small Rk} & \rotatebox{90}{\small Hv} \\
\midrule
\multirow{3}{*}{\rotatebox{90}{\resizebox{\ifdim\width>1.15cm 1.15cm\else\width\fi}{!}{\shortstack[c]{\small C.\ beam \\ \scriptsize $(\times 10^{-1})$}}}} & 12 & 4.88 & 5.36 & 5.29 & \underline{5.70} & \textbf{5.89} & 5.23 & 6.07 & 5.88 & \underline{6.21} & \textbf{6.33} & 5.20 & 6.32 & 6.13 & \underline{6.41} & \textbf{6.58} \\
 & 36 & 5.37 & 5.67 & 5.70 & \underline{5.91} & \textbf{5.93} & 5.60 & 6.22 & 6.18 & \underline{6.33} & \textbf{6.36} & 5.65 & 6.49 & 6.45 & \underline{6.51} & \textbf{6.61} \\
 & 60 & 5.56 & 5.77 & 5.79 & \underline{5.92} & \textbf{5.93} & 5.74 & 6.28 & 6.25 & \underline{6.33} & \textbf{6.36} & 5.83 & \underline{6.54} & 6.51 & 6.51 & \textbf{6.61} \\
\cmidrule{2-17}
\multirow{3}{*}{\rotatebox{90}{\resizebox{\ifdim\width>1.15cm 1.15cm\else\width\fi}{!}{\shortstack[c]{\small P.\ vessel \\ \scriptsize $(\times 10^{-1})$}}}} & 12 & 0.995 & 1.30 & 1.27 & \underline{1.57} & \textbf{1.94} & 1.83 & 2.31 & 2.21 & \underline{2.53} & \textbf{3.23} & 3.16 & 3.64 & 3.61 & \underline{4.17} & \textbf{4.97} \\
 & 36 & 0.946 & 1.74 & 1.65 & \underline{1.94} & \textbf{2.91} & 1.90 & 2.72 & 2.60 & \underline{3.13} & \textbf{4.50} & 3.34 & 4.31 & 4.20 & \underline{4.80} & \textbf{6.05} \\
 & 60 & 1.45 & 1.96 & \underline{2.15} & \underline{2.40} & \textbf{3.79} & 2.41 & 3.15 & 3.17 & \underline{3.62} & \textbf{5.04} & 3.82 & 4.82 & 4.56 & \underline{5.21} & \textbf{6.18} \\
\cmidrule{2-17}
\multirow{3}{*}{\rotatebox{90}{\resizebox{\ifdim\width>1.15cm 1.15cm\else\width\fi}{!}{\shortstack[c]{\small G02 \\ \scriptsize $(\times 10^{-4})$}}}} & 12 & 0 & 0 & 0 & 0 & 0 & 0 & 0 & 0 & 0 & 0 & 0 & 0 & 0 & 0 & 0 \\
 & 36 & 0 & 0 & 0 & 0 & 0 & 0 & 0 & 0 & 0 & 0 & 0 & 0 & 0 & 0 & 0 \\
 & 60 & 0 & 0 & 0 & 0 & 0 & 0 & 2.83 & 0 & 0 & 0 & 0 & 0 & 0 & 0 & 0 \\
\cmidrule{2-17}
\multirow{3}{*}{\rotatebox{90}{\resizebox{\ifdim\width>1.15cm 1.15cm\else\width\fi}{!}{\shortstack[c]{\small G06 \\ \scriptsize $(\times 10^{-1})$}}}} & 12 & 0.182 & \underline{\textbf{0.361}} & \underline{\textbf{0.534}} & \underline{\textbf{0.550}} & \textbf{0.470} & 0.421 & \underline{\textbf{1.34}} & \underline{\textbf{1.55}} & 1.17 & \textbf{1.43} & 0.862 & \underline{2.67} & \underline{2.62} & \underline{2.54} & \textbf{3.07} \\
 & 36 & 0.610 & \underline{1.62} & \underline{\textbf{1.86}} & \underline{1.55} & \textbf{2.04} & 1.45 & \underline{3.08} & \underline{\textbf{3.26}} & \underline{3.13} & \textbf{3.52} & 2.58 & \underline{4.54} & 4.25 & \underline{4.59} & \textbf{5.01} \\
 & 60 & 1.01 & 2.25 & \underline{\textbf{2.71}} & 2.34 & \textbf{2.74} & 2.26 & \underline{3.88} & \underline{3.81} & \underline{3.83} & \textbf{4.22} & 3.42 & \underline{5.05} & 4.59 & \underline{5.06} & \textbf{5.47} \\
\cmidrule{2-17}
\multirow{3}{*}{\rotatebox{90}{\resizebox{\ifdim\width>1.15cm 1.15cm\else\width\fi}{!}{\shortstack[c]{\small G08 \\ \scriptsize $(\times 10^{-1})$}}}} & 12 & 2.92 & 5.51 & 5.54 & \underline{5.56} & \textbf{5.59} & 3.05 & 7.10 & 7.13 & \underline{7.15} & \textbf{7.16} & 0.710 & 1.88 & \underline{1.98} & 1.87 & \textbf{2.00} \\
 & 36 & 3.07 & 5.57 & 5.58 & \underline{5.59} & \textbf{5.60} & 3.22 & 7.15 & 7.15 & \underline{7.16} & \textbf{7.17} & 0.750 & 1.98 & \underline{2.00} & 1.92 & \textbf{2.01} \\
 & 60 & 3.31 & 5.58 & 5.59 & \underline{5.59} & \textbf{5.60} & 3.57 & 7.16 & 7.16 & \underline{7.16} & \textbf{7.17} & 0.830 & 1.99 & \underline{2.01} & 1.93 & \textbf{2.01} \\
\cmidrule{2-17}
\multirow{3}{*}{\rotatebox{90}{\resizebox{\ifdim\width>1.15cm 1.15cm\else\width\fi}{!}{\shortstack[c]{\small G09 \\ \scriptsize $(\times 10^{-1})$}}}} & 12 & 4.37 & 4.75 & 4.76 & \underline{5.08} & \textbf{5.39} & 4.56 & 5.06 & 5.14 & \underline{5.43} & \textbf{5.70} & 4.52 & 5.22 & 5.30 & \underline{5.63} & \textbf{5.94} \\
 & 36 & 4.82 & 5.14 & 5.17 & \underline{5.41} & \textbf{5.59} & 4.94 & 5.45 & 5.44 & \underline{5.71} & \textbf{5.89} & 4.98 & 5.64 & 5.65 & \underline{5.92} & \textbf{6.11} \\
 & 60 & 4.97 & 5.26 & 5.26 & \underline{5.48} & \textbf{5.64} & 5.13 & 5.55 & 5.54 & \underline{5.77} & \textbf{5.93} & 5.19 & 5.74 & 5.74 & \underline{5.99} & \textbf{6.16} \\
\cmidrule{2-17}
\multirow{3}{*}{\rotatebox{90}{\resizebox{\ifdim\width>1.15cm 1.15cm\else\width\fi}{!}{\shortstack[c]{\small G24 \\ \scriptsize $(\times 10^{-1})$}}}} & 12 & 4.10 & 4.28 & 4.34 & \underline{4.57} & \textbf{4.65} & 4.29 & 4.56 & 4.66 & \underline{4.69} & \textbf{4.83} & 4.30 & 4.73 & 4.76 & \underline{4.87} & \textbf{4.97} \\
 & 36 & 4.51 & 4.72 & 4.75 & \underline{4.95} & \textbf{5.03} & 4.67 & 4.85 & 4.88 & \underline{5.02} & \textbf{5.08} & 4.68 & 4.94 & 4.94 & \underline{5.06} & \textbf{5.12} \\
 & 60 & 4.65 & 4.82 & 4.86 & \underline{5.03} & \textbf{5.10} & 4.79 & 4.93 & 4.96 & \underline{5.08} & \textbf{5.13} & 4.80 & 5.00 & 5.00 & \underline{5.11} & \textbf{5.15} \\
\midrule
\multicolumn{2}{r}{\small Wins} & \textbf{0} (0) & \textbf{1} (2) & \textbf{3} (4) & \textbf{1} (17) & \textbf{18} & \textbf{0} (0) & \textbf{1} (3) & \textbf{2} (3) & \textbf{0} (17) & \textbf{18} & \textbf{0} (0) & \textbf{0} (4) & \textbf{0} (4) & \textbf{0} (14) & \textbf{18} \\
\bottomrule
\end{tabular}
}
\end{table}
 
\FloatBarrier
\section{Discussion}
\label{sec:discussion}

The results on the fully constrained problem reinforce what other papers have already found, namely that preserving some infeasible solutions may improve, rather than hinder, search for the constrained optimum \citep{ray2009infeasibility,mezura2005simple}. Hence our claim that shadow-front information can essentially be obtained for free.

The results further allow a recommendation for each scenario of
Table~\ref{tab:scenarios}. Two
findings shape it. The hypervolume mode, although designed for joint
relaxation, is the most accurate mode, or tied for it, in most configurations
of every scenario. And it pays for this with a selection whose computational cost grows
steeply with the number of constraints and with the population size, and which
is by far the most expensive with three or four constraints and the largest population. The cheaper modes therefore
keep a role, which we describe for each scenario.

For the \emph{end points under single relaxation}, the hypervolume mode wins
most often at the low and medium levels and comes second to the relaxed mode,
by one configuration, at the high level. The joint relaxation box contains the
single-relaxation curves on its faces, and spreading the population over the
box also spreads it along those faces toward their far ends. Of the cheaper
modes, the relaxed mode, which was designed for this scenario, is the choice:\
it is never significantly worse than the baseline, it wins more often
than the ranking and joint modes, and it is the most accurate mode of all on
Mystery, together with the identical joint mode, and on G24 at the larger
population sizes. The ranking mode's design hypothesis,
that recovering the surrounding trade-off structure would also give a more
reliable estimate of the end points than independent round-robin picks,
holds only in part:\ the ranking mode is more accurate than the relaxed
mode on the cantilevered beam and G09, but less accurate on Branin, Mystery,
G24, and G08, mostly at the larger population sizes. With one constraint it has only two rankings, which are strongly
correlated, so there is little trade-off structure for non-dominated sorting to
exploit beyond what the relaxed mode already obtains directly.

For the \emph{whole curve under single relaxation}, the hypervolume mode is,
on our problems, again the best choice, and the ranking mode, which was
designed for this scenario, wins most often among the remaining modes at every level. It is also
the only mode that beats the hypervolume mode consistently on this metric,
on G24 in every configuration. The ranking mode's non-dominated sorting of rank
vectors spreads the population similarly to the hypervolume mode,
but with a weaker pressure toward joint relaxations. It keeps a further role:\ it
is the strongest optimizer of all modes on the best-so-far objective, and its
selection costs no more than non-dominated sorting with nearest-neighbor
niching. When coverage matters but the hypervolume mode's cost is too high,
the ranking mode is the choice.

For the two \emph{joint-relaxation} scenarios, the hypervolume mode is the
natural choice. It wins the hypervolume metric almost everywhere, which is a
property of its construction:\  it optimizes
the quantity being measured. More informative is that it is the most accurate
estimator of the threshold optimum, or tied for it, in every configuration of
the cantilevered beam, the pressure vessel, G08, and G09, and in six of G06.
On most of these problems the attainable gain grows
when several constraints are relaxed together, but the mode also leads where a
single constraint carries the joint gain, on the cantilevered beam at the low
level and on G08 at the low and medium levels. On the problems with several constraints, the relaxed mode, although it
relaxes one constraint at a time, estimates the joint end point better than the
hypervolume mode in only one configuration, of G02. The joint mode,
designed for this scenario, is the cheaper estimator of the joint end point:\
it comes second overall, ties with the hypervolume mode in four
configurations of the pressure vessel and two of G08, and is significantly more
accurate than it on G24 in three configurations at the larger population
sizes. The size of the joint gain does not decide either:\ on
G09 it exceeds the gain of $c_1$ by only about $5$ percent and the hypervolume
mode still leads clearly. A single run of the
hypervolume mode suffices for estimating the joint gain and, as the most
accurate mode, or nearly so, on the end points under single relaxation, also the per-constraint
gains. When the joint gain is well above the per-constraint gains, the
constraints interact and the joint surface is the object to study; when the two
agree, one constraint carries the gain.

The qualification, in every scenario, is computational cost. The relaxed, joint, and
ranking modes rank the pool by sorting and by non-dominated sorting, at a cost
comparable to that of the baseline. The hypervolume mode sorts the
threshold-feasible part of the pool into non-dominated fronts and orders each
front greedily, recomputing the hypervolume
contributions in $k+1$ dimensions after every removal, so its cost grows with
the size of the pool and steeply with the number of constraints. In our runs,
timed on one core of a laptop with cheap evaluations, a run of the relaxed or
joint mode takes at most about $1.8$ times as long as one of the baseline, a
run of the ranking mode at most about $2.2$ times, and a run of the
hypervolume mode at most about $1.4$ times on the problems with one or two
constraints. With three or four relaxed constraints, on the pressure vessel and
G09, a run of the hypervolume mode takes $1.9$ and $2.0$ times as long as one of
the baseline at $N = 36$, and $9.1$ and $6.0$ times as long at $N = 60$, about
$16$ and $9$ seconds against about $2$. This overhead amounts to up to about $1$
millisecond per
evaluation, negligible next to evaluations that take seconds or more, as in
simulation-based problems; but with several constraints it rises steeply with the population size, and with many constraints or cheap evaluations the
cheaper modes named for each scenario above are the better choice.
The hypervolume mode also spends selection pressure on spreading the population
over the relaxation box instead of pushing the strict optimum:\ on the
best-so-far objective it wins fewer configurations than the ranking and relaxed
modes at every level, and it is significantly worse than the baseline in every
configuration of G08, as are the joint and ranking modes, and in every
configuration of Branin with $N \geq 36$.

Across the scenarios, coverage and end-point accuracy are different things.
Hypervolume is best read as a diagnostic of how well a mode covers the
relaxed trade-off region, the quantity that matters in the whole-curve
scenarios, and the relative gap as the direct measure of end-point quality;
the ranking mode's position, second on coverage but, with the relaxed mode,
last of the alternative modes on the joint end points, shows that the two can diverge.

The gap metrics judge each mode against best-known references
(Section~\ref{subsec:exp_setup}), which are not certified values on nonconvex
problems. In practice the relaxation gain is read off a
single run as the difference between its best strictly feasible and its best
threshold-feasible value, and the two searches converge at different speeds.

The comparison against a feasibility-first baseline in the spirit
of~\citet{deb2000efficient} is a
comparison against conventional constraint-handling practice, not against a
competing estimator, because to our knowledge no prior method reports the
gain of relaxing a given constraint by a specified amount, or the
relaxation-versus-gain trade-off, as an output of a single run. The slope
estimate of~\citet{deb2010fast,deb2013biobjective} is no such estimator
either:\ it estimates the slope at zero violation of the front of objective
versus the sum of the normalized violations, a marginal rate used to set a
penalty parameter, not the gain of relaxing an individual constraint by a
specified amount. The closest alternative, re-optimizing per constraint level as in
\citet{saha2009ga}, requires one full run per level, whereas the shadow-front
here costs no evaluations beyond the run that finds the optimum. The other
obvious alternative, treating the constraints as additional objectives and
solving the multiobjective problem, does not do better either:\ our reference
front is such a run, with more than $150$ times the largest experiment budget,
and still falls short of the reference gain by up to $58$ percent at the joint
end points, while most runs of the modes beat it on G09 and on the pressure
vessel at the low and medium levels (Section~\ref{subsec:exp_setup}).

The shadow-front can also be read as a form of explanation of the optimization
outcome. It tells the decision maker which constraints bind at the optimum and
at what cost, which is the question sensitivity analysis has always answered
for linear and smooth convex problems. It serves what the framework for explainable AI
in operational research of~\citet{debock2024explainable} calls understandability:\
the ability of users to understand how a method works and the decisions it
reaches. In those terms, the shadow-front explains the
result in the vocabulary of the problem, its constraints and objective, rather
than the algorithm that produced it. In interactive multiobjective
optimization, explanations have so far concerned mainly the decision maker's
preferences and the trade-offs between the
objectives~\citep{misitano2022rximo,misitano2025emerging}. The survey
of~\citet{misitano2025emerging} distinguishes what-if explanations, of the
changes in a solution when, for example, the preferences or some decision
variables are modified, and names the explanation of further characteristics
of solutions, so far mainly their closeness to the constraint boundary, as an
open research question. The shadow-front is a what-if explanation of this kind
for the constraints:\ it shows how the attainable objective value changes when
a constraint is relaxed.

\section{Conclusions}
\label{sec:conclusions}

We adapted the shadow price to constrained single-objective black-box
optimization, where duality-based approaches are unavailable, by defining the
relaxation gain, the objective improvement obtained from a finite relaxation,
and the shadow-front, the trade-off between the amount by which a constraint is
relaxed, up to a tolerance that can be set by the decision maker,
and the attainable improvement in the objective. We considered four scenarios
with respect to what a decision maker might want to know about the benefit of
relaxing the constraints
up to this tolerance, according to whether only the end points or the
whole curve matter and whether constraints are relaxed one at a time or
jointly. We showed that in each scenario the information can be recovered as a
by-product of a generational evolutionary algorithm by modifying only its
environmental selection operator, with no additional objective evaluations and,
in most configurations, no loss in optimization performance.

We proposed four such operators, each designed for particular scenarios, and
compared them with a feasibility-first baseline on ten benchmark problems,
three population sizes, and three relaxation levels, with $100$ runs per
configuration. All four recovered the shadow-front and were significantly
worse than the baseline on the four scenario metrics in only one of $1224$
comparisons,
while matching or improving on the baseline's optimization performance in most
configurations. The hypervolume mode was, overall, the most accurate in all
four scenarios:\ it covered both the single-constraint curves and the joint
trade-off surface best, and it was the most accurate estimator of the end
points under single relaxation at the low and medium levels, second by one
configuration at the high level, and of the
threshold optimum on most problems. It is the choice in every scenario when its
selection cost, which grows steeply with the number of constraints, is
affordable. Of the cheaper modes, the relaxed mode was the better estimator of
the end points under single relaxation and the joint mode of the threshold
optimum, and each is the choice for its scenario, and the ranking mode, which needs no
hypervolume computation, was the strongest optimizer of all modes and the best
of the remaining modes at covering the trade-off region, and is the choice
when coverage matters.

\textbf{For a practitioner,} the recovered shadow-front turns constraint handling from a
verdict into decision-relevant information:\ from the same run and the same
evaluation budget, the decision maker receives, alongside a feasible solution,
an estimate of what relaxing each constraint by a given amount would buy. Our
method, experimental pipeline, and data are openly available, and the
experiments are fully reproducible.

Several limitations apply. We evaluate on ten textbook problems with cheap
evaluations and known optima; behavior on real-world problems remains to be
studied. The single-relaxation metrics
leave out the constraints that cannot improve the optimum on their own, so on
problems where only one constraint carries the gain the single and the joint
column measure much the same thing. We use DESDEO's standard NSGA-II operator settings and
do not tune them; with these settings some runs stall far from the optimum,
most often at the smaller population sizes, and on Townsend and G02 most runs
do, and other settings might change which mode
leads on individual problems. We have not compared the nearest-neighbor
niching of the ranking mode with other niching measures.
The hypervolume mode's per-generation cost is small next to the evaluations
when these are expensive, but it grows steeply with the number of constraints,
which limits the mode to problems with few constraints or expensive
evaluations.

Natural next steps are evaluation on a larger and more diverse problem set,
including simulator-coupled problems, in particular problems in which several
constraints carry a single-relaxation gain; different modes for different
constraints, for example the relaxed mode for some constraints and the joint
relaxation for others; adaptive choice
of niching in the ranking mode based on the dimension of the rank-vector space;
adaptive choice between the relaxed mode's per-constraint rankings and the
hypervolume mode's joint relaxation box, depending on whether one or several
constraints drive the attainable gain; and a theoretical characterization of
the per-mode bias in the empirical gain. Extending the approach to
multiobjective problems, where a decision maker would ask how a Pareto front or
a chosen solution changes under relaxation, is a further direction.

\section*{Data availability}

The run data, reference fronts, and derived results (tables, figures, and the
metrics behind them) are openly available on Zenodo~\citep{misitano2026data}.
The code implementing our proposed approach, and that produced the data and
regenerates the derived results, is archived on Zenodo~\citep{misitano2026code}.

\section*{CRediT authorship contribution statement}

\textbf{Giovanni Misitano:} Data curation, Formal analysis, Investigation, Methodology, Software,
Validation, Visualization, Writing -- original draft, Writing -- review \& editing.
\textbf{Juergen Branke:} Conceptualization, Investigation, Methodology, Validation, Writing --
original draft, Writing -- review \& editing.

\section*{Declaration of generative AI and AI-assisted technologies in the manuscript preparation process}

During the preparation of this work the authors used Claude Opus 5 and Claude
Opus 5.5 (Anthropic) to help refine the manuscript. All presented ideas and
analyses are the authors' own. After using these tools, the authors have always
reviewed and edited the content as needed. The authors stand behind every claim
made in the manuscript, and take full responsibility for the content of the
published article.

\section*{Acknowledgements}

Funding:\ This work was supported by the Research Council of Finland [grant number 355346] and
by the University of Jyv\"askyl\"a through a visiting fellowship of the second author. The funders
had no involvement in the study design, in the collection, analysis, and interpretation of data,
in the writing of the article, or in the decision to submit it for publication.

This work is part of the thematic research area Decision Analytics Utilizing Causal Models
and Multiobjective Optimization (DEMO, \url{https://www.jyu.fi/demo}) at the
University of Jyv\"askyl\"a.

\bibliographystyle{elsarticle-harv}
\bibliography{refs}

\end{document}